\input amstex
\documentstyle{amsppt}
\vsize=43pc
\magnification=\magstep 1
\def\picture #1 by #2 (#3){
    \vbox to #2{
        \hrule width #1 height 0pt depth 0pt
        \vfill
        \special{picture #3}}}

\def\scaledpicture #1 by #2(#3 scaled #4){{
   \dimen0=#1 \dimen1=#2
   \divide\dimen0 by 1000 \multiply\dimen0 by #4
   \divide\dimen1 by 1000 \multiply\dimen1 by #4
\picture \dimen0 by \dimen1 (#3 scaled #4)}}
\document
\topmatter
\title Complete Intersection  K-Theory and Chern Classes\endtitle
\author Satya Mandal\endauthor
\affil Institute of Mathematical Sciences
, C. I. T. Campus, Madras 600 113\endaffil
\email mandal@imsc.ernet.in\endemail
\dedicatory (Dedicated to my father)\enddedicatory
\date March 3, 1995\enddate
\endtopmatter
\def\refi #1 {\itemi{[#1]}}
\def\itemi{\par\hangi\textindenti}
\def\hangi{\hangindent=1.1cm}
\def\textindenti#1{\noindent\rlap{#1}\hskip1.1cm\ignorespaces}
\document
\head 0. Introduction\endhead
The purpose of this paper is to investigate the theory of complete intersection in
noetherian commutative rings from the K-Theory point of view.  (By {\it complete intersection
theory}, we mean questions like when/whether an ideal is the image of a projective module of
appropriate rank.)

The paper has two parts. In part one (Section 1-5), we deal with the relationship
between complete intersection and K-theory. The Part two (Section 6-8) is
, essentially, devoted to construction projectiove modules with certain cycles
as the total Chern class. Here Chern classes will take values in the
Associated graded ring of the Grothedieck $\gamma - filtration$ and as well
in the Chow group in the smooth case.

In this paper, all our rings are commutative and schemes are noetherian.
To avoid unnecessary complications, we shall assume that all our schemes
are connected.

For a noetherian scheme $X,~K_0(X)$ will denote the Grothendieck group
of locally free sheaves of finite rank over $X.$ Whenever it make sense,
for a coherent sheaf $M$ over $X,~[M]$  will denote the class
of $M$ in $K_0(X).$ We shall mostly be concerned with $X = Spec A$,
where $A$ is a noetherian commutative ring and in this case we shall
also use the notation $K_0(A)$ for $K_0(X).$    

\subhead Discussion on Part One (Section 1-5)\endsubhead
\bigskip
 
For a noetherian commutative ring $A$ of dimension $n$, we let 
$$F_0K_0A=\{[A/I]~ in~K_0A:~I~ is~a
~locally~complete~intersection~ideal~of~height~n\}.$$
In section 1, we shall prove that $F_0K_0A$ is a subgroup
of $K_0A$.  We shall call this subgroup $F_0K_0A$, the {\it zero cycle subgroup} of $K_0A$.  We
shall also see that (1.6), for a reduced affine algebra $A$ over an algebraically closed field
$k$, $F_0K_0A$ is the subgroup generated by smooth maximal ideals of height $n$.  The later
subgroup was considered by Levine \cite{Le} and Srinivas \cite{Sr}.

One of our main results (3.2) in Part One is that {\it for a noetherian commutative ring} $A$ {\it of
dimension}
$n$ {\it suppose that whenever} $I$ {\it is a locally complete intersection ideal of height} $n$
{\it with} $[A/I]=0$ {\it in} $K_0A$, {\it there is a projective (respectively} stably
free) $A$-module $P$ {\it of} rank $n$ {\it that maps onto} $I$. {\it Then for any
locally complete intersection ideal} $I$ {\it of} height $n$, {\it whenever}
$[A/I]$ {\it is divisible by}
$(n-1)!$ {\it in}
$F_0K_0A$,
$I$ {\it is image of a projective}
$A$-module
$Q$ {\it of} rank $n$ {\it (respectively} with $(n-1)!([Q]-[A^n]) = -[A/I]$ in $K_0A$).

In \cite{Mu2,(3.3)}, Murthy proved that for a reduced affine algebra $A$ over an algebraically closed
field $k$, for an ideal $I$, if $I/I^2$ is generated by $n=\ $dim $A$ elements then $I$ is image of
a projective $A$-module of rank $n$.  

In example (3.6), we show that for the coordinate ring
$$A=\Bbb R [X_0,X_1,X_2,X_3]/(X_0^2+X^2_1+X^2_2+X_3^2-1)$$ 
of real 3-sphere, the ideal
$I=(X_0-1,X_1,X_2,X_3)A$ is not the image of a projective A-module of rank 3, although $[A/I]=0$
in $K_0A$.

Another interesting result (3.4) in this part is that {\it suppose that}
$f_1,f_2,\dotsc,f_r$ {\it is a regular sequence in a noetherian commutative
ring} $A$ {\it of dimension} $n$ {\it and let} $Q$ {\it be a projective}
$A$-module {\it of rank} $r$ {\it that maps onto} $f_1,\dotsc,f_{r-1},
f_r^{(r-1)!}.$ {\it Then} $[Q] = [Q_0\oplus A]$ {\it in} $K_0A$ {\it
for some projective} $A$-module $Q_0$ {\it of rank} $r-1.$

When $dim~A~=~n~=~r~=~rank~Q$, this result(3.4)  has interesting comparison with
 the corresponding theorem of Mohan Kumar \cite{Mk1} for 
reduced affine algebras $A$ over algebraically closed fields.

More generally we prove that (3.5) {\it suppose} $A$ {\it is a noetherian
commutative ring of} dimension $n$ {\it and let} $J$ {\it be a locally complete intersection
ideal of} height
$r\leq n$.  {\it Assume that} $K_0A$ {\it has no} $(r-1)!$ {\it torsion,} $[A/J]=0$ {\it and}
$J/J^2$ {\it has free generators of the form} $f_1,f_2,\dotsc,f_{r-1},f^{(r-1)!}_r$ {\it in}
$J$.  {\it Let}
$Q$ {\it be a projective}
$A$-module {\it of} rank $r$ {\it that maps onto} $J$.
  {\it Then} $[Q]=[Q_0\oplus A]$ {\it in}
$K_0A$, {\it for some projective} $A$-module
$Q_0$ {\it of} rank $r-1$.

For reduced affine algebras $A$ of dimension $n$ over
algebraically closed fields $k$, and for $n=r$,
(3.5) is a consequence of the theorem of Murthy 
\cite{Mu2, Theorem 3.7}.  Besides these results
\cite{Mk1, Mu2}  (3.4) and (3.5)
are the best in this context , even for affine 
algebras over algebraically closed fields.
In fact, there is almost no result available in the case when
rank is strictly less than the dimension of the ring.

In section 1, we define and describe the zero cycle subgroup $F_0K_0A$ of
$K_0(A)$.  In section 2, for $k' =\Bbb Z$ or a field, we define the ring 
$$A_n=A_n(k')=\frac
{k'[S,T,U,V,X_1,\dotsc ,X_n,Y_1,\dotsc ,Y_n]}{(SU+TV-1,X_1Y_1+\cdots +X_nY_n-ST)}.$$
For our
purposes, $A_n$ serves like a "universal ring."  Besides doing the construction of the
"universal projective module" (2.6), we compute the $K_0A_n$, the Chow Group of $A_n$ and we
comment on the higher $K$-groups of $A_n$.

All the results in  Section 3 discussed above follows from a key Theorem (3.1).
In Section 4, we give the proof of (3.1).  In
Section 5, we give some more applications of (3.1).

\subhead Discussion on Part Two (Section 6-8)\endsubhead  
\bigskip

The purpose of this part of the the paper is to construct  
projective modules of appropriate 
rank that have certain cycles as its Chern classes and to 
consider related questions.

 For a noetherian scheme $X$
 of dimension $n$, $\Gamma(X)=\bigoplus_{i=1}^{n}\Gamma^{i}(X)$ will denote
 the the graded ring associated to the Grothendieck $\gamma-filtration$ of the
 Grothendieck group $K_{0}(X)$ and $CH(X) = \bigoplus_{i=1}^{n}CH^{i}(X)$ 
will denote the Chow group of cycles of  $X$ modulo rational equivalence.

\medskip

Our main construction (8.3) is as follows: {\it suppose} $X = Spec A$ {\it is a
 Cohen-Macaulay scheme of dimension} $n$ {\it and} $r \geq r_{0}$ 
{\it are  integers  with} $2r_{0}\geq n$ and $n \geq r .$
{\it Given a projective} $A$-module $ Q_{0}$ {\it of rank} $r_{0}-1$ 
{\it and a sequence of locally
 complete intersection ideals} $I_{k}$ {\it of height} $k$ {\it for}
 $k = r_{0}$ {\it to} $r$ {\it such that}

(1) {\it the restriction} $Q_{0}|Y$ {\it is trivial for all locally complete intersection
 subschemes} $Y$ {\it of codimension at least} $r_{0}$ {\it and} 

(2){\it for} $k = r_{0}$ {\it to} $r~   I_{k}/I_{k}^{2}$ {\it has 
a free set of generators of the type}
\newline $f_{1},\ldots,f_{k-1},f_{k}^{(k-1)!}$
{it in} $I_{k}$, 

{\it then there is a projective} $A$-module $Q_{r}$ {\it of rank} $r$ {\it such that} 

(1) {\it for} $1\leq k \leq r_{0}-1$ {\it the} $kth$ {\it Chern class of}
 $Q_{0}$ {\it and} $Q_{r}$
{\it  are same and} 

(2){\it for} $k$ {\it between} $r_{0}$ {\it and} $r$ 
{\it the} $kth$  {\it Chern class of} $Q_{r}$ 
{\it is given by the cycle of} $A/I_{k}$, {\it upto a sign}.(Here Chern classes take values
in $\Gamma(X)\bigotimes{\Cal Q}$ or in the Chow group, if $X$ is nonsingular over
a field).{\it More precisely, we have}
$$[Q_{r}]-r = [Q_{0}]-(r_{0}-1)+\sum_{k=r_{0}}^{r}[A/J_{k}]$$ 
{\it in}
 $K_{0}(X)$ {\it where} 
$J_{k}$ {\it is a locally complete intersection ideal of height}
 $k$ {\it with} $[A/I_{k}]=
-(k-1)![A/J_{k}]$.

\medskip

Inductive arguments are used to do the construction of $Q_{r}$ in theorem(8.3).
Conversely, we prove theorem(8.2):

\medskip

{\it let} $A$ {\it be a commutative noetherian ring of
dimension} $n$ {\it and} $X=Spec A$.{\it Let} $J$ {\it be a locally complete intersection ideal of 
height} $r$ {\it so that} $J/J^{2}$ {\it has a free set of generators of the form} 
$f_{1},f_{2}, \ldots ,f_{r-1},f_{r}^{(r-1)!}$. {\it Let} $Q$ 
{\it be projective} $A$-module  {\it of 
rank} $r$ {\it that maps onto} $J$.{\it Then there is a projective} $A$
-module $Q_{0}$ {\it of rank}
$r-1$ {\it such that the first} $r-1$ {\it Chern classes of}
 $Q$ {\it and} $Q_{0}$ {\it are same.}(Here
again Chern classes take values in $\Gamma (X)$  or
in the Chow group,if $X$ is nonsingular
over a field). {\it Infact, if} $K_{0}(X)$  {\it is torsion free then}
$[Q_{0}]$ {\it is unique in} $K_{0}(X)$.

\medskip

Both in the statements of theorem (8.2) and (8.3), we considered locally
complete intersection ideals $J$ of height $r$ so that $J/J^{2}$ has free
set of generators of the form $f_{1},f_{2},\ldots ,f_{r-1},f_{r}^{(r-1)!}$ in
$J$. For such an ideal $J,~ [A/J]=(r-1)![A/J_{0}]$, for some locally complete
intersection ideal $J_{0}$. Consideration of such ideals are supported in
theorem(8.1):

\medskip

{\it Let} $A$ {\it be a noetherian commutative ring of dimension} $n$ {\it and}
$X=Spec A$. {\it Assume that} $K_{0}(X)$ {\it has no} $(n-1)!$-torsion.
{\it Let} $I$ {\it be 
locally complete intersection ideal of height} $n$ {\it that is image of a projective}
$A$-module $Q$ {\it of rank} $n$.{\it Also suppose that} $Q_{0}$ {\it is an}
 $A$-module {\it of rank}
$n-1$ {\it so that the first} $n-1$ {\it Chern classes of}
 $Q$ {\it and} $Q_{0}$ {\it are same.Then}
$[A/I]$ {\it is divisible by} $(n-1)!$.

\medskip

For a variety $X$,what cycles of $X$,in $\Gamma(X)$ or in the Chow group,that may
 appear as the total Chern class of a locally free sheaf of appropriate rank
had always been an interesting question, although not much is known in this
direction.

\medskip

For affine smooth three folds $X=Spec A$ over algebraically closed fields,
Mohan Kumar and Murthy\cite{MM} proved that(see 8.9) if $c_{k}$ is a
cycle in $CH^{k}(X)$, for $k=1,\,2,\,3$ then there are projective $A$-modules
 $Q_{k}$ of rank
$k$ so that 
\newline $(1)total~Chern~ class~ of~ Q_{1}\, is\, 1+c_{1}$
\newline $(2) the~ total~ Chern~ class~ of~ Q_{2}~ is~ 1+c_{1}+c_{2}$,
\newline $(3) the~ total~ Chern~ class~ of~ Q_{3}~ is~ 1 +c_{1}+c_{2}+c_{3}$.
\newline

We give a stronger version (8.10) of this theorem(8.9) of 
Mohan Kumar and Murthy \cite{MM}. Our theorem (8.10) applies
to any smooth three fold $X$ over any field such that $CH^3(X)
$ is divisible by 2.

\medskip

Murthy\cite{Mu2} also proved that if $X=SpecA$ is a smooth affine variety 
of dimension $n$ over an algebraically closed filed $k$ and 
$c_{n}$ is a codimension $n$ cycle in the Chow group of $X$
then there is a projective $A$-module $Q$ of rank $n$ so that the total
Chern class of $Q$ is $1+c_{n}$. 
We give a stronger version(8.7) of this theorem of 
Murthy \cite{Mu2}. This version(8.7) of the theorem applies to all smooth affine
varieties $X$ of dimension $n$, over any field, so that $CH^{n}(X)$ 
is divisible by $(n-1)!$.

\medskip

Murthy \cite{Mu2} also proved : suppose that $X=SpecA$ is a smooth affine
variety of dimension $n$ over an algebraically closed field $k$. For 
$i\, =\, 1\,$ to $n$ let $c_{i}$ be a codimension $i$ cycle in the Chow group of $X$. Then
there is a projective $A$-module $Q_{0}$ of rank $n-1$ with total Chern class
$1+c_{1}+\cdots+c_{n-1}$ if and only if there is a projective $A$-module $Q$ of
rank $n$ with total Chern class $1+c_{1}+\cdots +c_{n}$.
We also give an alternative proof (8.8)
of this theorem of Murthy \cite{Mu2} .

\medskip

Besides these results \cite{MM,Mu2} not much else is known in this direction.
Our results in section 8 apply to any smooth affine variety over any field
and also consider codimesion $r$ cycles where $r$ is strictly less than
the dimension of the variety. Consideration of Chern classes in the 
Associated graded ring of the Grothendieck $\gamma-filtration$ in
the nonsmooth case is, possibly, the only natural thing to do because
the theory of Chern classes in the Chow group is not available in 
such generality. Such consideration of Chern classes in the Associated
graded ring of the Grothendieck $\gamma-filtration$ was never done  
before in this area .

In section 6, we set up the notations and other formalism about the
Grothendieck Gamma filtration, Chow groups and Chern classes.
In this section  we also give an
example of a smooth affine variety $X$ for which the Grothedieck Gamma filtration
of $K_{0}(X)$ and the filtration by the codimension of the support do not agree.

In section 7, we set up some more preliminaries. Our main results 
of the Part Two of the paper are in section 8. 
\medskip

I would like to thank M. P. Murthy for the innumerable number of discussions I had with him over a
long period of time.  My sincere thanks to M. V. Nori for many stimulating discussions.  I would
also like to thank Sankar Dutta for similar reasons. I thank D. S.
Nagaraj for helping me to improve the exposition and for many discussions.

\head Part One : Section 1-5\endhead
\head Complete Intersection and K-Theory\endhead

\medskip

In this part, we investigate the relationship between complete intersction
and K-theory.

\medskip

\head 1. The Zero Cycle Subgroup\endhead
For a noetherian commutative ring $A$, $K_0(A)$ will denote the Grothendieck group of projective
$A$-modules of finite rank.  We define $F_0K_0A=\{[A/I]$ in $K_0A:I$ is a locally complete
intersection ideal of height $n$ = dim$A$\}.

In this section, we shall prove that $F_0K_0A$ is a subgroup of $K_0A$.  We call this subgroup
$F_0K_0A$, the {\it zero cycle subgroup}.  We shall also prove that if $A$ is an affine algebra
over an algebraically closed field $k$, this notation $F_0K_0A$ is consistent with the notation
used by Levine \cite{Le} and Srinivas \cite{Sr} for the subgroup of $K_0A$ generated by $[A/\Cal
M]$, where $\Cal M$ is a smooth maximal ideal in $A$.

\proclaim{Theorem 1.1}Suppose $A$ is a noetherian commutative ring of dimension  $n$.  Then
$F_0K_0A$ is a subgroup of $K_0A$.\endproclaim

The proof of (1.1) will follow from the following Lemmas.

\proclaim{Lemma 1.2} $F_0K_0A=-F_0K_0A$.\endproclaim

\demo{Proof} Suppose $I$ is a locally complete intersection ideal of height $n$ = dim $A$ and
$x=[A/I]$ is in $F_0K_0A$.

Let $I=(f_1,f_2,\dots ,f_n)+I^2$.  By induction, we shall find $f'_1,f'_2,\dots
,f'_r$ in $I$ such that 

(1) $(f'_r,\dots ,f'_r,f_{r+1},\dots ,f_{n})+I^2=I$,

(2) $(f'_1,\dots ,f'_r)$ is a regular sequence.

Suppose we have picked $f'_1,f'_2,\dots, f'_r$ as above and $r<n$.  Let $\frak
p _1, \dots, \frak p_s$ be the associated primes of $(f'_1,f'_2,\dots
,f'_r)$. If $I$ is contained in $\frak p _1$, then height $\frak p _1=n$ and since $I_{{\frak
p}_1}$ is complete intersection of height $n$, $A_{{\frak p}_1}$ is Cohen-Macaulay ring of height
$n$.  This contradicts that $\frak p _1$ is associated prime of
$(f'_1,\dots ,f'_r)$. So, $I$ is not contained in $\frak p _i$ for $i=1$ to $s$. Let
$\{P_1,\dots ,P_t\}$ are maximal among $\{\frak p _1,\dots ,\frak p _s\}$ and assume that $f_{r+1}$
is in $P_1,\dots ,P_{t_0}$ and not in $P_{t_0+1},\dots ,P_t$. Let $a$ be in $I^2 \cap P_{t_0+1}
\cap\dots\cap P_t\setminus P_1\cup P_2\cup\dots\cup P_{t_0}$.  Let $f'_{r+1}=f_{r+1}+a$. Then
$f'_{r+1}$ does not belong to $P_i$ for $i=1$ to $t$ and hence also does not belong to $\frak
p _1,\dots ,\frak p _s$. Hence we have that

(1) $(f'_1,\dots ,f'_r,f'_{r+1},f_{r+2},\dots ,f_n)+I^2=I$ and

(2) $f'_1,\dots ,f'_r,f'_{r+1}$ is a regular sequence.

Therefore, we can find a regular sequence $f'_1,f'_2,\dots,f'_n$ such that
$I=(f'_1,\dots ,f'_n)\allowmathbreak +I^2$. So, $(f'_1,f'_2,\dots
,f'_n)=I\cap J$ for some ideal
$J$ with $I+J=A$. Since $f'_1,\dots ,f'_n$ is a regular sequence, $J$ is locally complete
intersection ideal and $[A/I]+[A/J]=[A/(f'_1,f'_2,\dots ,f'_n)]=0$. Hence $[A/J]=-x$ is
in $F_0K_0A$. So, the proof of (1.2) is complete.\enddemo

\proclaim{Lemma 1.3} Suppose $A$ is a noetherian commutative ring of height $n$ and $I$ is a
locally complete intersection ideal of height $n$.  Let $\Cal M_1,\dots ,\Cal M_k$ be maximal
ideals that does not contain $I$. There are $f_1,f_2,\dots ,f_n$ such that 

(1) $f_1,\dots ,f_n$ is a regular sequence,

(2) $I=(f_1,\dots ,f_n)+I^2$,

(3) for a maximal ideal $\Cal M$, if $(f_1,\dots ,f_n)$ is contained in $\Cal M$, then $\Cal M
\ne \Cal M _i$ for $i=1$ to $k$.\endproclaim

\demo{Proof} As in the proof of (1.2), we can find a regular sequence $f_1,\dots ,f_n$ such that
$I=(f_1,\dots ,f_n)+I^2$. We readjust $f_n$ to avoid $\Cal M_1,\dots ,\Cal M_k$ as follows. Let
$\frak p_1,\dots ,\frak p_s$ be the associated primes of $(f_1,\dots ,f_{n-1})$. Then $I$ is not
contained in $\frak p_i$ for $i=1$ to $s$. Let $\{P_1,P_2,\dots ,P_t\}$ be maximal among $\{\frak
p_1,\dots ,\frak p_s,\Cal M_1,\dots ,\Cal M_k\}$. Assume that $f_n$ is in $P_1,\dots ,P_{t_0}$
and not in $P_{t_0+1},\dots ,P_t$. Let $a$ be in $I^2 \cap P_{t_0+1}\cap\dots\cap P_t\setminus
P_1\cup P_2\cup\dots\cup P_{t_0}$ and $f'_n=f_n+a$. Then $f'_n$ is not in $\Cal M_1,\dots
,\Cal M_k$. So,

(1) $f_1,f_2,\dots ,f_{n-1}, f'_n$ is a regular sequence

(2) $I=(f_1,f_2,\dots ,f_{n-1},f'_n)+I^2$ and

(3) if $(f_1,\dots ,f_{n-1}, f'_n) \subseteq \Cal M$ for a maximal ideal $\Cal M$ then $\Cal M
\ne \Cal M_i$ for $i=1$ to $k$.

This completes the proof of (1.3)\enddemo

\proclaim{Lemma 1.4} Let $A$ be as in \rom{(1.1)}. Then $F_0K_0A$ is closed under
addition.\endproclaim

\demo{Proof} Let $x$ and $y$ be in $F_0K_0A$. Then $x=[A/I]$ and by (1.2), $y=-[A/J],$ where $I$
and $J$ are locally complete intersection ideals of height $n$. Let $\{\Cal M_1,\dots ,\Cal M_k\}
= V(I)\setminus V(J)$, the maximal ideals that contain $I$ and do not contain $J$. By (1.3),
there is a regular sequence $f_1,f_2,\dots ,f_n$ such that $J=(f_1,\dots ,f_n)+J^2$ and for
maximal ideals $\Cal M$ that contains $(f_1,\dots ,f_n), \Cal M \ne \Cal M_i$ for $i=1$ to
$k$.

Let $(f_1,\dots ,f_n)=J\cap J'$, where $J+J'=A$ and $J'$ is a locally complete
intersection ideal of height $n$. Then $y=-[A/J]=[A/J']$. Also note that $I+J'=A$. Hence
$x+y=[A/I]+[A/J']=[A/IJ']$, and $IJ'$ is a locally complete intersection ideal of height
$n$. So the proof of (1.4) is complete.\enddemo

Clearly, the proof of Theorem (1.1) is complete by (1.2) and (1.4). Now we proceed to prove that
for a reduced affine algebras $A$ over a field $k$, $F_0K_0A$ is generated by regular points.

\proclaim{Theorem 1.5} Suppose $A$ is a reduced affine algebra over a field $k$ of dimension $n$.
Then $F_0K_0A$ is generated by the classes $[A/\Cal M ]$, where $\Cal M$ runs through all the
regular maximal ideals of height $n$.\endproclaim

\demo{Proof} Since the regular locus of $A$ is open (see [K]), there is an ideal $J$ of $A$ such
that $V(J)$ is the set of all prime ideals $P$ such that $A_P$ is not regular. Since $A$ is
reduced, height $J \geq 1$.

Let $G$ be the subgroup of $K_0A$, generated by all classes $[A/\Cal M ]$, where $\Cal M$ is a
regular maximal ideal of $A$ of height $n$. Clearly $G$ is contained in $F_0K_0A$. Now let
$x=[A/I]$ be in $F_0K_0A$, with $I$ a locally complete intersection ideal of height $n$. Let
$I=(f_1,f_2,\dots ,f_n)+I^2$. By induction we shall find $f'_1,f'_2,\dots ,f'_r$ for $r \leq n$,
such that 

(1) $I=(f'_1,\dots ,f'_4,f_{r+1},\dots ,f_n)+I^2$,

(2) $(f',_1,f'_2,\dots ,f'_r)$ is a regular sequence and 

(3) for a prime ideal $P$ of $A$, if $J+(f'_1,\dots ,f'_r)$ is contained in $P$ then either
height $P>r$ or $I$ is contained in $P$.

We only need to show the inductive step.  Suppose we have picked $f'_1,\dots ,f'_r$ as
above. Let $\frak p_1, \dots \frak p_k$ be the associated primes of $(f'_1,\dots f'_r)$
and let $Q_1, Q_2, \dots ,Q_s$ be the minimal primes over $(f'_1, f'_2, \dots ,f'_r)+J$
so that $I$ is not contained in $Q_i$ for $i=1$ to $s$. So, we have height $Q_i>r$. As before, we
see that $I$ is not contained in $\frak p_i$ for $i=1$ to $k$.

Let $\{P_1,P_2,\dots ,P_t\}$ be the maximal elements in $\{\frak p_1,\dots ,\frak p_k,Q_1,\dots
,Q_s\}$ and let $f_{r+1}$ be in $P_1, \dots ,P_{t_0}$ and not in $P_{t_0+1},\dots ,P_t$. Let $a$
be in $I^2 \cap P_{t_0+1}\cap\dots\cap P_t\setminus P_1\cup P_2\cup\dots\cup P_{t_0}$. Write
$f'_{r+1}=f_{r+1}+a$. Then $f'_{r+1}$ will satisfy the requirement.

Hence we have sequence $f'_1,f'_2,\dots ,f'_n$ such that 

(1) $I=(f'_1,\dots ,f'_n)+I^2$,

(2) $f'_1,f'_2,\dots ,f'_n$ is a regular sequence and 

(3) if a maximal ideal $\Cal M$ contains $(f'_1,\dots f'_n)+J$ then $I$ is contained in
$\Cal M$. If $(f'_1,\dots ,f'_n)=I\cap I'$, then $I+I'=A$ and $I'$ is a locally complete
intersection ideal of height $n$. Also, if a maximal ideal $\Cal M$ contains $I'$, then
$\Cal M$ is a regular maximal ideal of height $n$. So, $[A/I']$ is in $G$ and hence
$x=[A/I]=-[A/I']$ is also in $G$. The proof of (1.5) is complete.\enddemo

\remark{Remark 1.6} From the proof of (1.5), it follows that (1.5) is valid for any noetherian
commutative ring $A$ such that the singular locus of spec $A$ is contained in a closed set $V(J)$
of codimension at least one. Similar arguments works for smooth ideals.\endremark

\head{2. The Universal Constructions}\endhead

For $k'=\Bbb Z$ or a field, we
let$$K=K(k^{')}=\frac{k'[S,T,U,V]}{(SU+TV-1)}$$$$A_n=A_n(k')=\frac{k'[S,T,U,V,X_1,\dots
,X_n,Y_1,\dots ,Y_n]}{(SU+TV-1,X_1Y_1+\dots
+X_nY_n-ST)}$$$$B_n=B_n(k')=\frac{k'[T,X_1,\dots ,X_n,Y_1,\dots ,Y_n]}{(X_1Y_1+X_2Y_2+\dots
+X_nY_n-T(1+T))}$$

\flushpar By the natural map $A_n\rightarrow B_n$, we mean the map that sends $T\rightarrow T,
S\rightarrow 1 + T, U \rightarrow 1, V\rightarrow -1$. (We will continue to denote the images of
upper case letter variables in $A_n$ or $B_n$, by the same symbol).

The ring $B_n$, was considered by Jouanlou \cite{J}. Later $B_n$ was further used by Mohan Kumar
and Nori \cite{Mk2} and Murthy \cite{Mu2}. The purpose of this section is to establish that
$A_n$ behaves much like $B_n$.

\subhead The Grothendieck Group and the Chow Group of $A_n$\endsubhead

\medskip

For a ring $A$ and $X$ = Spec $A, K_0(A)$ or $K_0(X)$ ({\it respectively} $G_0(A)$ or $G_0(X)$)
will denote the Grothendieck Group of finitely generated projective modules ({\it respectively}
finitely generated modules) over $A$. $CH^k(A)$ or $CH^k(X)$ will denote the Chow Group of cycles
of codimension $k$ modulo rational equivalence and $CH(X)=\bigoplus CH^k(X)$ will denote the
total Chow group of $X$.

\proclaim{Proposition 2.1} Let $\lambda _n=[A_n/(X_1,X_2,\dots ,X_n,T)]$ in $G_0(A_n)$. Then
$G_0(A_n)$ is freely generated by $\epsilon _n=[A_n]$ and $\lambda _n$. In fact, the natural map
$G_0(A_n)\rightarrow G_0(B_n)$ is an isomorphism.\endproclaim

\demo{Proof} We proceed by induction on $n$. If $n=0$, then $A_0\approx A_0/(S)\times A_0/(T)$.
Since $A_0/(S)\approx A_0/(T)\approx k'[S^{\pm 1},V]$, the proposition holds in this case.

Now assume $n>0$. We have $A_{n_{X_n}}\approx K[X_1,\dots ,X_{n-1},X^{\pm 1}_n,Y_1,\dots
,Y_{n-1}]$ and $A_n/(X_n)\approx A_{n-1}[Y_n]$. Since $G_0(K)\approx \Bbb Z$ (see \cite{Sw1}, \S
10),
$G_0(A_{n_{X_n}})\approx \Bbb Z$ and also by induction $G_0(A_n/(X_n))\approx
G_0(A_{n-1}[Y_n])\approx G_0(A_{n-1})$ is generated by $[A_{n-1}]$ and $\lambda _{n-1}$. Now we
have the exact sequence $G_0(A_n/(X_n))$ $@> i_* >>$ $G_0(A_n)$ $@> j^* >>$
$G_0(A_{n_{X_n}})\rightarrow 0$. Since the $i_*(\lambda_{n-1})=\lambda_n$ and $i_*([A_{n-1}])=0,
G_0(A_n)$ is generated by $\lambda_n$ and $[A_n]$.

It is also easy to see that the natural map $G_0(A_n)\rightarrow G_0(B_n)$ sends $\lambda_n$ to
$\beta_n=[B_n/(X_1,\dots ,X_n,T)]$ and $[A_n]$ to $[B_n]$. Since $\beta_n$ and $[A_n]$ are free
generators of $G_0(B_n)$ (see \cite{Sw1, \S 10}/\cite{Mu2}), $\lambda_n$ and $[A_n]$ are free
generators of
$G_0(A_n)$. This completes the proof of (2.1).\enddemo

\proclaim{Proposition 2.2} Let $\lambda '_n$ be the cycle defined by $A_n/(X_1,\dots X_n,T)$
in $CH(A_n)$ and let $\epsilon '_n=[Spec\ A_n]$ be the cycle of codimension zero.  Then
$CH(A_n)$ is freely generated by $\epsilon '_n$ and $\lambda'_n$. That means $CH^j(A_n)=0$
for $j\ne 0, n, CH^0(A_n)=\Bbb Z \epsilon'_n\approx\Bbb Z$ and $CH^n(A_n)=\Bbb Z
\epsilon '_n\approx
\Bbb Z$.\endproclaim

Before we prove (2.2), we prove the following easy lemma.

\proclaim{Lemma 2.3} $CH(K) = \Bbb Z [Spec\ K]$.\endproclaim

\demo{Proof} Note that $K/UK\approx k'[T^{\pm 1},S]$ and $K_U\approx k'[T,U^{\pm 1},V]$. Now
the lemma follows from the exact sequence $CH_j(K/UK)\rightarrow CH_j(K)\rightarrow
CH_j(K_U)\rightarrow 0$ for all
$j$. (Here we use the notation $CH_j(S)=CH^{dim\ x-j}(X))$.\enddemo

\demo{Proof of (2.2)} For $n=0, A_0=A_n\approx K/(ST)\approx K/(S)\times K/(T)\approx
k'[T^{\pm 1},U]\times k'[S^{\pm 1},V].$ $CH^0(A_0)\approx \Bbb Z[V(S)]\oplus\Bbb Z
[V(T)]\approx \Bbb Z\epsilon'_0\oplus\Bbb Z\lambda'_0$ and $CH^j(A_0)=0$ for all $i>0$.
So, the Proposition 2.2 holds for $n=0$. Let $$d_n=dim\ A_n=\cases 2n+2\ \text{if}\ 
k'\text{is a field}\\2n+3\ \text{if}\ k'=\Bbb Z .\endcases$$

Let $n>0$ and assume that the proposition holds for $n-1$. Since $A_{n_{X_n}}\approx
K[X_1,\dots,X_{n-1},X^{\pm 1}_n,Y_1,\dots ,Y_{n-1}]$ and $A_n/(X_n)\approx A_{n-1}[Y_n],
CH(A_{n_{X_n}})\approx \mathbreak\Bbb Z [Spec\ A_{n_{X_n}}]$ and $CH(A_n/(X_n)\approx
CH(A_{n-1}[Y_n])\approx CH(A_{n-1})$. By induction it follows that $CH^j(A_n/(X_n)=0$ for
$j\ne 0, n-1$ and $CH^{n-1}(A_n/(X_n)$ is freely generated by $[A_n/(X_1,\dots ,X_n,T)]$.

Now consider the exact sequence $CH^{j-1}(A_n/(X_n))\rightarrow
CH^j(A_n)\rightarrow CH^j(A_{n_{X_n}})\allowmathbreak\rightarrow 0$.  It follows that for $j\ne
0, n, CH^j(A_n)=0$ and clearly, $CH^0(A_n)$ is freely generated by $\epsilon'_n=[Spec\ A_n]$.
Also
$CH^n(A_n)$ is generated by the image of \linebreak $[A_n/(X_1,\dots ,X_n,T)]$, which is
$\lambda '_n$. Since the natural map $CH^n(A_n)\rightarrow G_0(A_n)$ maps $\lambda'_n$ to
$\lambda_n$ and
$\lambda_n$ is a free generator, it follows that $\lambda'_n$ is also torsion free.  Hence
$CH^n(A_n)=\Bbb Z \lambda'_n\approx \Bbb Z$. This completes the proof of (2.2).\enddemo

\subhead Higher $K$-Groups of $A_n$\endsubhead

\medskip

Much of this section is inspired by the arguments of Murthy \cite{Mu3} and Swan \cite{SW1}. 
Again for a ring $A$, $G_i(A)$ will denote the $i$-th $K$-group of the category of finitely
generated $A$-modules.  For a subring $K$ of $A$, $\widetilde G_i(A,K)$ will denote the cokernel
of the map $G_i(K)\rightarrow G_i(A)$.  As also explained in \cite {SW1}, if there is an
augmentation $A\rightarrow K$ with finite tordimension, then $0\rightarrow G_i(K)\rightarrow
G_i(A)\rightarrow 
\widetilde G_i(A,K)\rightarrow 0$ is a split exact sequence.

Following is a remark about the higher $K$-groups of $A_n$.

\proclaim{Theorem 2.4} Let $k$ be a field or $\Bbb Z$ and $K=K(k)$ and $A_n=A_n(k)$. Then 

(1) $G_i(K)\approx G_{i}(k)\oplus G_{i-1}(k)$ for all $i\ge 0$

(2) for $i\ge 0, n\ge 1,\widetilde{G}_i(A_{n,k}) \approx \widetilde{G}_i(A_1,K)$ and
$0\rightarrow G_i(K)\rightarrow G_i(A_n)\rightarrow \widetilde{G}_i(A_n,K)\rightarrow 0$ is
split exact,

(3) There is a long exact sequence $\dots \rightarrow G_i(k[T^{\pm 1}])\oplus G_i(k[S^{\pm
1}])\rightarrow G_i(A_1)\rightarrow G_i(K[X^{\pm 1}])\overset\partial\to\longrightarrow 
G_{i-1}(k[T^{\pm 1}])\oplus G_{i-1}(k[S^{\pm
1}])\rightarrow G_{i-1}(A_1))\rightarrow\dots$\endproclaim

\demo{Proof} The statement (1) is a theorem of Jouanolou \cite{J}. To prove (2), note that
all rings we consider are regular and that $K\rightarrow A_n$ has an augmentation for
$n\ge 1$. Also note that for $n\ge 2, K[X_n]\rightarrow A_n$ is a flat extension. So, it
induces a map of the localization sequences

$$\matrix&&&G_i(K)&&&&&\\&&&ss&&&&&\\
\dots\rightarrow\!\!&G_i(K)&\!\!\rightarrow\!\!&G_i(K[X_n])&\!\!\rightarrow\!\!&G_i(K[X_n^{\pm
1}])&\rightarrow&G_{i-1}(K)&\rightarrow\dots\\
&\downarrow&&\downarrow&&\downarrow ss&&\downarrow&\\
\dots\rightarrow\!\!&G_i(A_n/X_n)&\!\!\rightarrow\!\!&G_i(A_n)&\!\!\rightarrow\!\!&
G_i(A_{n_{X_n}})&\!\!\rightarrow\!\!&G_i(A_n/X_n)&\!\!\rightarrow\dots\\ &ss&&&&&&ss&\\
&G_i(A_{n-1})&&&&&&G_{i-1}(A_{n-1})&\endmatrix\tag2.5$$

\flushpar Also note that $0\rightarrow G_i(K)\rightarrow G_i(K[X^{\pm 1}])\rightarrow
G_{i-1}(K)\rightarrow 0$ is a split exact sequence $([Q])$. This will induce an exact
sequence $\rightarrow G_i(K)\rightarrow G_i(A_{n-1})\rightarrow
\widetilde{G}_i(A_n,K)\rightarrow G_{i-1}(K)\rightarrow\dots$. Since $G_i(K)\rightarrow
G_i(A_{n-1})$ is split exact, it follows that $0\rightarrow G_i(K)\rightarrow
G_i(A_{n-1})\rightarrow\widetilde{G}_i(A_n,K)\rightarrow 0$ is split exact.  Hence
$\widetilde{G}_i(A_{n-1},K)\approx \widetilde{G}_i(A_n,K)$. This establishes statement (2).
The statement (3) is the is immediate consequence of the
 localization sequence $([Q])\rightarrow
G_i(A_1/(X_1)\rightarrow G_i(A_1)\rightarrow G_i(A_{1_{X_1}})\rightarrow
G_{i-1}(A_1/(X_1))\rightarrow\dots$ .  This completes the proof of (2.4).\enddemo

\subhead Construction of the Universal Projective Module\endsubhead

\medskip

Under this subheading we construct a projective module over $A_n$, which will be useful
in the later sections. This construction is similar to the construction of Mohan Kumar and
M. V. Nori \cite {Mk2} over $B_n$.

\proclaim{Proposition 2.6} Let $J_n$ be the ideal $(X_1,X_2,\dots ,X_n,T)$ in $A_n$. Then

(1) for $n\ge 3$, there is no projective $A_n$-module of rank $n$ that maps onto $J_n$.

(2) There is a projective $A_n$-module $P$ of rank $n$ that maps onto the ideal
$J'_n=(X_1,\dots,X_{n-1})A_n+J^{(n-1)!}_n$ such that $([P])-n)=[A/J_n]=-\lambda_n$ in
$G_0(A_n)$.\endproclaim

\demo{Proof} The proof of the statement (1) is similar to the argument in \cite {Mk2} or
this can also be seen by tensoring with $B_n$ and using the result in \cite {Mk2}.

To prove statement (2), note that $J_{n_{S}}=(X_1,\dots ,X_n)$ and
$J'_{n_{S}}=(X_1,\dots ,X_{n-1},\allowmathbreak X^{(n-1)!}_n)$. Also, since $(X_1,\dots
,X_{n-1},X^{(n-1)!}_n)$ is an unimodular row in $A_{n_{ST}}$, by Suslin's Theorem
(\cite {S}), there is an $n\times n$-matrix $\gamma$ in $\Bbb M_n(A_n)$ such that $det(\gamma
)=(ST)^u$ for some $u\ge 0$ and the first column of $\gamma$ is the transpose of
$(X_1,X_2,\dots ,\allowmathbreak X_{n-1},X^{(n-1)!}_n)$. Let $f_1:A^n_{n_{S}}\rightarrow
J'_{n_{S}}$ be the map that sends the standard basis  $e_1,\dots ,e_n$ of $A^n_{n_{S}}$ to
$X_1,X_2,\dots ,X_{n-1},X^{(n-1)!}_n$ and let $f_2:A^n_{n_{T}}\rightarrow
J_{n_{T}}\approx A_{n_{T}}$ be the map that sends the standard basis $e_1,\dots ,e_n$ to
$1,0,0,\dots ,0$. As in the paper of Boratynski \cite {B}, by patching $f_1$ and $f_2$ by
$\gamma$, we get a surjective map $P\rightarrow J'_n$, where $P$ is a projective
$A_n$-module of rank
$n$.

Now we wish to establish that $([P]-[A^n_n])=-\lambda_n$. By tensoring with $\Bbb Q$, in
case $k=\Bbb Z$, we can assume that $k$ is a field.

The rest of the argument is as in Murthy's paper \cite {Mu2}.  Let $[P]-[A^n_n]=m\lambda
n$. So, $C_n([P]-[A^n_n])=(-1)^n[V(J'_n)]=(-1)^n(n-1)!\lambda'_n$. Also, by the
Riemann-Roch theorem, $C_n([P]-[A^n_n])=m
C_n(\lambda_n)=m(-1)^{n-1}(n-1)!\lambda'_n$. Hence it follows from (2.2) that $m=-1$.
Hence ${P}-[A^n_n]=-\lambda_n$.\enddemo

\head 3. The Main Results in Part One\endhead

Our main results follows from the following central theorem.

\proclaim{Theorem 3.1} Let $A$ be a commutative noetherian ring of dimension $n$ and let
$I$ and $J_0$ be two ideals that contain nonzero divisors and $I+J_0=A$. Assume that
$J_0$ is a locally complete intersection ideal of height $r$ with
$J_0=(f_1,\dots,f_r)+J^2_0$ and let $J=(f_1,\dots,f_{r-1})+J^{(r-1)!}_0$. Suppose $Q$ is
a projective $A$-module of rank $r$ and $\varphi :Q\rightarrow IJ$ is a surjective map. 
Then

(i) there is a projective $A$-module $P$ of rank $r$ that maps onto $J$ with
$[P]-[A^r]=-[A/J_0]$in$K_0(A)$;

(ii) further, there is  a surjective map from $Q\oplus A^r$ onto $I\oplus P$;

(iii) in particular, there is a projective $A$-module $Q'$ of rank $r$ that maps onto
$I$ and $[Q']=[Q]+[A/J_0]$ in $K_0(A)$.\endproclaim

In the rest of this section we shall use this theorem (3.1) to derive its main
consequences and the proof of (3.1) will be given in the next section.

\proclaim{Theorem 3.2} Let $A$ be a noetherian commutative ring of dimension $n\ge 1$.
Also assume that for locally complete intersection ideals $I$ of height $n$, whenever
$[A/I]=0$ in $K_0(A)$, $I$ is an image of a projective (\rom{respectively}, with stably
free) $A$-module $Q$ of rank $n$.

Then for locally complete intersection ideals $I$ of height $n$, if $[A/I]$ is divisible
by $(n-1)!$ in $F_0K_0A$ then $I$ is image of a projective $A$-module
$Q'$ (\rom{respectively}, with $(n-1)!([Q']-n)=-[A/I])$ of rank $n$.\endproclaim

\remark{Remark 3.3} If $A$ is a reduced affine algebra over an algebraically closed field,
Murthy \cite{Mu2} proved that for any ideal $I$ of $A$, if $I/I^2$ is generated by
$n=$ dim $A$ elements then $I$ is an image of a projective $A$-module of rank $n$.\endremark

\demo{Proof of \rom{(3.2)}} Let $I$ be a locally complete intersection ideal of height $n$, so
that $[A/I]$ is divisible by $(n-1)!$ in $F_0K_0A$. Let $[A/I]=(n-1)![A/J]$ in $F_0K_0A$.

Let $\Cal M_1,\dots,\Cal M_k$ be maximal ideals that contains $I$ and that does not contain
$J$. By Lemma (1.3), we can find a locally complete intersection ideal $J'$ of height $n$
such that $[A/J]=-[A/J']$ and $I+J'=A$. Now let $J'=(f_1,...,f_n)+J^{'2}$ and
$J^{''}=(f_1,\dots,f_{n-1})+J^{'(n-1)!}$. So, $[A/I]=(n-1)![A/J]=-(n-1)![A/J^{"}]=-[A/J^{''}]$
and
$I+J^{''}=A$. Hence $[A/IJ^{''}]=0$. By hypothesis,there is a projective ({\it respectively},
stably free) $A$-module $Q$ of rank $n$ that maps onto $IJ^{''}$. By (3.1) there is a
projective module $Q'$ of rank $r$ that maps onto $I$ and $[Q']=[Q]+[A/J']$. Hence also
$(n-1)!([Q']-[Q])=(n-1)![A/J']=-[A/I]$. So,the proof of (3.2) is complete.\enddemo

Our next two applications (3.4, 3.5) of (3.1) are about splitting projective modules.

\proclaim{Theorem 3.4} Let $A$ be noetherian commutative ring and let
$f_1, f_2,\ldots , f_r$ be a regular sequence. Let $Q$ be projective $A-module$
of rank $r$ that maps onto $(f_1,\ldots ,f_{r-1},f_r^{(r-1)!})$. Then
$[Q] = [Q_0\oplus A]$ for some projective $A-module~Q_0$ of rank $r-1$. 
\endproclaim

The proof of (3.4) is immediate from (3.1) by taking $J_0 = (f_1,f_2,
\ldots ,f_r)$ and $I = A.$

Following is a more general version of (3.4).

\proclaim{Theorem 3.5} Let $A$ be a noetherian commutative ring of dimension $n$ and let $J$
be a locally complete intersection ideal of height $r\ge 1$, such that $J/J^2$ has free
generators of the type $f_1,f_2,\dots,f_{r-1},f^{(r-1)!}_r$ in $J$. Suppose $[A/J]=0$ in
$K_0(A)$ and assume $K_0(A)$ has no $(r-1)!$ torsion. Then, for a projective $A$-module $Q$ of
rank $r$, if $Q$ maps onto $J$, then $[Q]=[Q_0]+1$ in $K_0(A)$ for some projective $A$-module
$Q_0$ of rank $r-1$.\endproclaim

\demo{Proof} First note that we can assume that $f_1, f_2, \ldots ,f_r$is
a regular sequence. 
We can find an element $s$ in $J$ such that $s(1 + s)=f_1g_1+f_2g_2+\dots
+f_{r-1}+f^{(r-1)!}_rg_r$ for some $g_1,g_2,\dots g_r$. Let $J_0=(f_1,f_2,\dots
,f_{r-1},f_r,s)$. Then $J_0$ is a locally complete intersection ideal of height $r$ and
$J=(f_1,\dots,f_{r-1})+J^{(r-1)!}_0$. Let $I=A$. Then by (3.1), there is a projective
$A$-module $Q'$ of rank $r$ that maps onto $A$ and $[Q']=[Q]+[A/J_0]=[Q]$. Since $Q'=Q_0\oplus
A$ for some $Q_0$, the theorem (3.5) is established.\enddemo

\remark{Remark 3.6} For reduced affine algebras $A$ over algebraically closed fields $k$,
Murthy \cite{Mu2} proved a similar theorem for $r=n=$ dim $A\ge 2$. In that case, if chark
$=0$ or chark $=p\ge n$ or $A$ is regular in codimension $1$, then $F_0K_0A$ has no
$(n-1)!$-torsion. (See \cite{Le}, \cite{Sr}, \cite{Mu2})\endremark

Before we close this section we give some examples.

\example{Example 3.7} Let $A=\Bbb R [X_0,X_1,X_2,X_3]/X^2_0+X^2_1+X^2_2+X^2_3-1)$ be the
coordinate ring of the real $3$-sphere $S^3$. Then $K_0(A)=\Bbb Z$ (see \cite{Hu} and
\cite{Sw2}) and $CH^3(A)=\Bbb Z /2\Bbb Z$ generated by a point \cite{CF}. Since, in this case
for any projective $A$-module $Q$ of rank $3$, the top Chern Class $C_3(Q)=0$ in $CH^3(A)$, no
projective $A$-module will map onto the ideal $I=(X_0-1,X_1,X_2,X_3)A$. This is a situation,
when $[A/I]=0$ in $K_0(A)$, but $I$ is not an image of a projective module of rank
$3$.\endexample

\head 4. Proof of Theorem 3.1\endhead

In this section we give the proof of Theorem 3.1.  First we state the following easy lemma.

\proclaim{Lemma 4.1} Suppose $A$ is a noetherian commutative ring and $I$ and $J$ are two
ideals that contain nonzero divisors.  Let $I+J=A$. Then we can find a nonzero divisor $s$ in
$I$ such that $(s,J)=a$.\endproclaim

Now we are ready to prove (3.1).

\demo{Proof of \rom{(3.1)}} The  first part of the proof is to find a nonzero divisor $s$ in $I$
such that 

(1) $(s,J_0)=A$,

(2) after possibly modifying $f_1,\dots ,f_r$, we have $sJ_0\subseteq (f_1,\dots ,f_r)$ and

(3) $Q_s$ is free with basis $e_1,\dots ,e_r$ such that $\varphi_s(e_1)=f_1,\dots
,\varphi_s(e_{r-1})=f_{r-1}$ and $\varphi_s(e_r)=f^{(r-1)!}_r$.

First note that there is a nonzero divisor $s_1$ in $I$ such that $As_1+J=A$. Now let $\Cal
P_1,\Cal P_2,\dots ,\Cal P_r$ be the associated primes of $A_{s_1}$ such that $\Cal
P_i+J_{s1}=A_{s_1}$. We pick maximal ideals $\Cal M_1,\dots ,\Cal M_k$ in spec $(A_{s_1})$ such
that $\Cal P_i\subseteq\Cal M_i$ for $i=1$ to $k$ and let $\Cal M_0=\Cal M_1\cap\dots\cap\Cal
M_k$. Then $J_{0_{s_1}}+\Cal M_0=A_{s_1}$. Let $a+b=1$ for $a$ in $J^2_{s1}$ and $b$ in $\Cal
M_0$. Let $f'_r=bf_r+a$. It follows that
$J_{0_{s_1}}=(f_1,\dots,f_{r-1},f'_r)+J^2_{0_{s_1}}\Cal M$.

Hence there is $s_2=1+t_2$ in $1+J_{0_{s_1}}\Cal M_0$, such that
$J_{0_{s_1s_2}}=(f_1,\dots,f'_r)$. Clearly, $s_2$ is not in $\Cal P_1,\Cal P_2,\dots ,\Cal P_r$.
If $\Cal P$ is any other associated prime of $A_{s_1}$ and $s_2=1+t_2$ is in $\Cal P$, then
$J_{0_{s_1}}+\Cal P =A_{s_1}$, which is impossible. So, we have found a nonzero divisor $s_2$
in $1+J_{0_{s_1}}\Cal M_0$ such that $J_{0_{s_1s_2}}=(f_1,f_2,\dots ,f'_r)$.

Now let $K$ be the kernel of $\varphi_{s_1s_2}:Q_{s_1s_2}\rightarrow J_{s_1s_2}$. Since
$J_{s_1s_2}=(f_1,f_2,\dots ,f_{r-1},\allowmathbreak f^{'(r-1)!}_r)$, there are $e'_1,e'_2,\dots
,e'_r$ in $Q_{s_1s_2}$ such that $\varphi (e'_1)=f_1,\dots ,\varphi
(e'_{r-1})=f_{r-1},\allowmathbreak
\varphi (e'_r)=f^{'(r-1)!}_r$.

By tensoring $0\rightarrow K\rightarrow Q_{s_1s_2}\rightarrow J_{s_1s_2}\rightarrow 0$ by
$A_{s_1s_2}/\Cal M_{0s_2}$, we get an exact sequence $0\rightarrow K/\Cal M_0K\rightarrow
Q_{s_1s_2}/\Cal M_0Q_{s_1s_2}\overset\bar\varphi\to\rightarrow J_{s_1s_2}/J_{s_1s_2} \Cal
M_0\approx A_{s_1s_2}/\Cal M_0{s_2}\rightarrow 0$
and $\overline{\varphi (e'_r)}=\overline{f^{(r-1)!}}$ is a unit in $A_{s_1s_2}/\Cal
M_{0_{s_2}}$. (Bar means module $\Cal M_{0_{s_2}}$). So there are $E_1,E_2,\dots ,E_{r-1}$ in
$K$, such that images of $E_1,E_2,\dots ,E_{r-1},e'_r$ is a basis of $Q_{s_1s_2}/\Cal M_0
Q_{s_1s_2}$.

Write $e_1=be'_1+aE_1,e_2=be'_2+aE_2,\dots ,e_{n-1}=be'_{r-1}+aE_{r-1},e_r=e'_r$. It is easy
to see that $e_1,\dots ,e_r$ is a basis of $Q_{s_1s_2W}$, where $W=1+J_{s_1s_2}\Cal M_0$. So,
there is $s_3=1+t_3$ in $1+J_{s_1s_2}\Cal M_0$ such that $e_1,\dots ,e_r$ is a basis of
$Q_{s_1s_2s_3}$. As before, $s_3$ is a nonzero divisor in $A_{s_1s_2}$. Of course,
$\varphi_{s_1s_2s_3}(e_1)=bf_1,\dots ,\varphi_{s_1s_2s_3}(e_{r-1})=bf_{r-1},
\varphi_{s_1s_2s_3}(e_r)=(f'_r)^{(r-1)!}$. By further inverting a nonzero divisor in
$1+J_{0s_1s_2}$, we can also assume that $bf_1,\dots ,bf_{r-1},f'_r$ generate $J_{0s_1s_2s_3}$.

So, we are able to find a nonzero divisor $s$ in $A$ and a free basis $e_1,e_2,\dots ,e_r$ of
$Q_s$ such that, after replacing $f_1$ by $bf_1,\dots f_{r-1}$ by $bf_{r-1}$ and $f_r$ by
$f'_r$, we have 

(1) $s$ is in $I$ and $su+t=1$ for some $t$ in $J_0$ and $u$ in $A$.

(2) $sJ_0\subseteq (f_1,\dots ,f_r)$

(3) $\varphi (e_1)=f_1,\dots\varphi (e_{r-1})=f_{r-1},\varphi (e_r)=f^{(r-1)!}_r$.

We had to go through all these technicalities because we wanted to have a nonzero divisor $s$.
Now let $st=g_1f_1+g_2f_2+\dots +g_rf_r$ and let $s^kQ\subseteq\overset r\to{\underset i=1\to
\bigoplus} Ae_i\approx A^r$ for some $k\ge 0$.

By replacing $Q$ by $s^kQ$ and $I$ by $s^kI$, we can assume that

(4) $s^{k+1}$ is in $I$,

(5) $ts=g_1f_1+\dots +g_rf_r$.

(6) There is an inclusion $i:Q\rightarrow A^r=Ae_1+\dots +Ae_r$ such that $Q_s=A^r_s$ and

(7) $\varphi_s(e_i)=f_i$ for $i=1$ to $r-1$ and $\varphi_s(e_r)=f^{(r-1)!}_r$.

Let $A_r=A_r(\Bbb Z)$ be as in section 2 and let us consider the map $A_r\rightarrow A$ that
sends $X_i$ to $f_i$, $Y_i$ to $g_i$ for $i=1$ to $r$ and $T$ to $t$, $S$ to $s$, $U$ to $u$
and $V$ to $1$. By the theorem of Suslin (\cite{S}) there is an $r\times r$ matrix $\gamma$ in
$\Bbb M_r(A_r)$ with its first column equal to the transpose of $(X_1,X_2,\dots ,X_{r-1},
X^{(r-1)!}_r)$ and with det$(\gamma)=(ST)^a$ in $A_r$, for some integer 
$a\ge 1$. Now let $\alpha$ be
the image of $\gamma$ in $\Bbb M_r(A)$.

We shall consider $\alpha$ as a map $\alpha :A^r\rightarrow A^r$ and let
$\alpha_0:Q\rightarrow A^r$ be the restriction of $\alpha$ to $Q$.

Define the $A$-linear map $\varphi_0:A^r\rightarrow A$ such that $\varphi_0(e_i)=f_i$ for
$i=1$to $r-1$ and $\varphi_0(e_r)=f^{(r-1)!}_r$. Also let $\varphi_1=(1,0,\dots
,0):A^r_t\rightarrow A_t$ be the map defined by $\varphi_1(e_1)=1$ and $\varphi_1(e_i)=0$ for
$i=2$ to $r$.
Also let $\varphi_2=(\varphi_0)_s$.
Note that $\varphi :Q\rightarrow IJ$ is the restriction of $\varphi_0$ to $Q$ and hence the
diagram
$$\CD Q_t@>\varphi>>IJ_t\\
@VV\alpha_0 V        @VVV \\
A^r_t @>\varphi_1>> A_t\endCD$$

is commutative.

Now consider the following fibre product diagram:

\vskip7pt
\centerline{\scaledpicture 5.24in by 4.35in (picture1 scaled 600)}
\vskip6pt

Here $\varphi_2=(\varphi_0)_s$ is a surjective and the maps $\eta$ and $\psi$ on the upper left
hand corner are given by the properties of fibre product diagram.

Clearly, the map $\psi:P\rightarrow J$ is surjective.  Further, since $\alpha$ is the image of
$\gamma$ it follows from (2.6) that $[P]-[A^r]=$ image of $-\lambda_r=-[A/J_0]$.

Now it remains to show that $Q\oplus A^r$ maps onto $I\oplus P$.

Note that the diagram
$$\CD Q @>\phi>> IJ \\
@VV\eta V  @VVV \\
P @>>> J\endCD$$

\flushpar is commutative because it is so on $D(t)$ and $D(s)$. Also note that the map
$\eta_s:Q_s\rightarrow P_s$ is an isomorphism and hence $s^pP$ is contained in $\eta(Q)$ for
some $p\ge 1$.

Write $K=$ kernel$(\psi)$. So, the sequence $0\rightarrow K\rightarrow
P\overset\psi\to\rightarrow J\rightarrow 0$ is exact.

Since $Tor_i(J,A/s^pA)=0$, the sequence $0\rightarrow K/s^pK\rightarrow P/s^pP\rightarrow
J/s^pJ\approx A/s^pA\rightarrow 0$ is exact. In the following commutative diagram of exact
sequences$$\CD P/s^pP @>\bar{\psi}>> J/s^pJ\approx A/s^pA\rightarrow 0\\@VVV @VVV \\
A^r_t/s^pA^r_t @>\bar{\varphi}_1>> J_t/s^pJ_t\rightarrow 0,\endCD$$

\flushpar the vertical maps are isomorphism. But since $\bar{\varphi}_1=(1,0,\dots ,0),
K/s^pK=$ kernel $\bar{\psi}\approx$ ker $\bar{\varphi}_1$ is a free $A/s^pA$-module of rank $r-1$.

Now write $M=$ kernel $\phi$. Then we have the following commutative diagram

\newpage
$$0\rightarrow K\rightarrow P\overset\psi\to\longrightarrow J\rightarrow 0$$
$$\uparrow\quad\quad\uparrow\eta\quad\quad\uparrow$$
$$0\rightarrow M\rightarrow Q\longrightarrow IJ\rightarrow 0$$

\flushpar of exact sequences.

Define the map $\delta : P\oplus I\rightarrow J+I=A\rightarrow 0$ such that $\delta(p,x)=\psi
(p)-x$ for $p$ in $P$ and $x$ in $I$ and let $L=$ kernel$(\delta)$. So, $0\rightarrow
L\rightarrow P\oplus I\rightarrow J+I=A\rightarrow 0$ is an exact sequence and $L\oplus A$ is
isomorphic to $P\oplus I$. Th, it is enough to show that $Q\oplus A^{r-1}$ maps onto
$L$.

But $L$ is isomorphic to $\psi^{-1}(IJ)$ and we have the following commutative diagram of
exact sequences:$$0\rightarrow K\rightarrow L\overset\psi\to\longrightarrow J\rightarrow 0$$
$$\uparrow\quad\quad\uparrow\eta\quad\quad\vert\vert$$
$$0\rightarrow M\rightarrow Q\longrightarrow IJ\rightarrow 0$$

Note that $s^pK$ is contained in $\eta (M)$. So $K/s^pK$ maps onto $K/\eta (M)$. Therefore
$K/\eta (M)$ is generated by $r-1$ elements.

As $Q\oplus K/\eta (M)$ maps onto $L, Q\oplus A^{r-1}$ maps onto $L$. This completes the proof
of Theorem 3.1.\enddemo

\head 5. The Theorem of Murthy\endhead

In this section we give some applications of (3.1), which was inspired by the fact that the
Picard group of smooth curves over algebraically closed fields are divisible.

\proclaim{Theorem 5.1} Let $A$ be a commutative ring of dimension $n$ and $I$ be a locally
complete intersection ideal of height $n$ in $A$. Suppose that $I$ contains a locally complete
intersection ideal $J'$ of height $n-1$ and there is a projective $A$-module $Q_0$ of rank
$n-1$ that maps onto $J'$. If the image of $I$ in $A/J'$ is invertible and is divisible by
$(n-1)!$ in \rom{Pic}$(A/J')$, then there is a projective $A$-module $Q$ of rank $n$ that
maps onto $I$ and $(n-1)!([Q]-[Q_0]-[A])=-[A/I]$ in $K_0(A)$.\endproclaim

\demo {Proof} Let bar \lq\lq--" denote images in $A/J'$. Since $\bar{I}$ is divisible by
$(n-1)!$ in Pic$(A/J')$, it's inverse is also divisible by $(n-1)!$. Let $J$ be an ideal of $A$
such that $J'\subseteq J, I+J=A$ and $\bar{J}^{(n-1)!}=\bar{I}^{-1}$ in Pic$(A/J')$. Hence
$\overline{IJ}^{(n-1)!}=(j',f)/J'$ for some $f$ in $I$.

Write $G=J'+J^{(n-1)!}$. We can also find a $g$ in $J$ such that $J=(J',g)+J^2$. Since
$J'/J'J$ is locally generated by $(n-1)$ elements, $J'/J'J$ is $(n-1)$-generated. Let
$g_1,\dots ,g_{n-1}$ generate $J'/J'J$. We can find an element $s$ in $J$, such that
$J'_{1+s}=(g_1,\dots ,g_{n-1})$ and $J_{1+s}=(g_1,\dots ,g_{n-1}, g)$. Hence it also follows
that $G_{1+s}=(g_1,\dots ,g_{n-1},g^{(n-1)!})$. Therefore $G=(g_1,\dots ,g_{n-1})+J^{(n-1)!}$
and $J=(g_1,\dots ,g_{n-1},g)+J^2$. Since $\overline{IJ}^{(n-1)!}=(J',f)/J'$, it follows that
$IG=(f,J')$. As $Q_0\oplus A$ maps onto $IG$, by Theorem 3.1, there is a surjective map
$\varphi :Q_0\oplus A^{n+1}\rightarrow I\oplus P$, where $P$ is a projective $A$-module of
rank $n$ with $[P]-[A^n]=-[A/J]$. Let $Q=\varphi^{-1}(I)$. Then $Q$ maps onto $I$ and $Q\oplus
P\approx Q_0\oplus A^{n+1}$. So, $[Q]-[Q_0]-[A]=-([P]-[A^n])=[A/J]$. Hence
$(n-1)!([Q]-[Q_0]-[A])=(n-1)![A/J]=[A/G]=-[A/I]$. This completes the proof of (5.1)\enddemo

\proclaim{Corollary 5.2} Let $A$ be a smooth affine domain over an infinite field $k$ and let
$X=$ \rom{Spec} $A$. Assume that for all smooth curves $C$ in $X$, \rom Pic $C$ is divisible by
$(n-1)!$. If $I$ is a smooth ideal of \rom{height} $n =$ \rom{dim} $X$, then there is a
projective module $Q$ of \rom{rank} $n$, such that $Q$ maps onto $I$ and
$(n-1)!([Q]-[A^n])=-[A/I]$.\endproclaim

\demo{Proof} We can find elements $f_1,\dots f_{n-1}$ in $I$ such that $C=$ Spec$(A/(f_1,\dots
,f_{n-1}))$ is smooth \cite {Mu2, Corollary 2.4}. Now we can apply (5.1) with $Q_0=A^{n-1}$.\enddemo

\remark{Remark} Unless $k$ is an algebraically closed field, 
there is no known example of affine smooth variety that
satisfy the hypothesis of (5.2) about
the divisibility of the Picard groups. When $k$ is 
an algebraically closed field, (5.2) is a theorem of
Murthy \cite{Mu2, Theorem 3.3}
.\endremark

\head Part Two : Section 6-8\endhead
\head Projective Modules and Chern Classes\endhead

\medskip

This part of the paper is devoted to construct Projective modules with
certain cycles as the total Chern class and to consider related questions.  
Our main results in this Part are in section 8.

\bigskip

\head 6.  Grothendieck $\gamma-filtration$ and Chern Class formalism
\endhead

As mentioned in the introduction, for a noetherian scheme $X$, $K_{0}(X)$ will
denote the Grothendieck group of locally free sheaves of finite rank over $X$.
All schemes we consider are connected and has an ample line bundle on it.

In this section we shall recall some of the formalisms about the Gorthendieck
$\gamma-filtrations$ of the Grothendieck groups and about Chern classes.The
main sources of this material are \cite {SGA6},\cite{Mn} and \cite{FL}.

\medskip

\example{Definitions and Notations 6.1} Let $X$ be noetherian scheme of dimension
$n$ and let $K_{0}(X)[[t]]$ be the power series ring over $K_{0}(X)$.Then

\medskip

a)$\lambda_{t} = 1+t\lambda^{1}+t^{2}\lambda^{2} +\cdots $ will denote the
additive to multiplicative group homomorphism from $K_{0}(X)$ to $1+tK_{0}(X)[[t]]$
induced by the exterior powers, that is $\lambda ^{i}([E]) = [\Lambda ^{i}(E)]$
for any locally free sheaf $E$ of finite rank  over $X$, and $i= 0,1,2,\ldots$,

\medskip

b)$\gamma_{t} = 1+t\gamma^{1}+t^{2}\gamma^{2} +\cdots$ will denote the map
$\lambda_{t/1-t}$, which is also an additive to multiplicative group homomorphism.

\medskip

c)We let $F^{0}K_{0}(X)= K_{0}(X)$,$F^{1}K_{0}(X)= Kernel(\epsilon)$ where
$\epsilon :K_{0}(X) \rightarrow {\Bbb Z}$ is the rank map.

For positive integer $k, F^{k}K_{0}(X)$ will denote the subgroup
of $K_{0}(X)$ generated by the elements
$\gamma ^{k_{1}} (x_{1}) \gamma ^{k_{2}}(x_{2}) \ldots \gamma ^{k_{r}}(x_{r})$ such 
that $\sum_{i=1}^{r} k_{i} \geq k $ and $x_{i}$ in $F^{1}K_{0}(X)$.
We shall often write $F^{i}(X)$ for $F^{i}K_{0}(X)$.

Recall that
$$F^{0}(X) \supseteq F^{1}(X) \supseteq F^{2}(X) \supseteq \cdots $$ 
is the Grothendieck $\gamma$-filtration of $K_{0}(X)$. Also note that $F^{n+1}(X)=0$ (see \cite {FL,Mn}).

\medskip

d)$\Gamma(X) = \bigoplus_{i=0}^{n}\Gamma ^{i}(X)$ will denote the graded
ring associated to the Grothendieck $\gamma-filtration$.
\newline If $x$ is in $F^{k}(X)$, then the image of x in $\Gamma^{k}(X)$ will
be called the $cycle~of~x$ and be denoted by $Cycle(x)$.

\medskip

e) For a locally free sheaf $E$ of rank $r$ over $X$ and for nonnegative
interger $i$, the $ith$ Chern class of $E$ is defined as
$$c_{i}(E) = \gamma^{i}([E]-r) ~~modulo ~~F^{i+1}(X).$$

This will induce a Chern class homomorphism
$$c_{t}:K_{0}(X) \rightarrow 1+ \bigoplus_{i=1}^{n}\Gamma^{i}(X)t^{i}$$
which is also an additive to multiplicative group homomorphism. We write
$$c_{t}(x) = 1+\sum_{i=1}^{n}c_{i}(x)t^{i}$$ with $c_{i}(x)$ in $\Gamma^{i}(X)$.

\medskip

f) We recall some of the properties of this Chern class homomorphism:
\newline
(1) if $x$ is in $F^{k}(X)$, then 
$$c_{i}(x) = 0~~  for ~~	1\leq i < r$$
$$c_{r}(x) = (-1)^{r-1}(r-1)! Cycle(x),$$
\newline
(2) if $E$ is a locally free sheaf of rank $r$ and if $E$ maps onto a locally complete
intersection sheaf $I$ of ideals of height $r$ then $[{\Cal O}_{X}/I] = \sum_{i=0}^{r}(-1)^{i} \lambda^{i}[E]$
is in $F^{r}(X)$ and 
$$c_{r}([E]) = (-1)^{r}Cycle ([{\Cal O}_{X}/I]) $$

\endexample

	Now we shall set up some notations about the formalism of Chern classes 
in Chow groups.

\example{Notations~and~Facts 6.2} Let $X$ be a noetherian scheme of dimension $n$ and let
$$CH(X) = {\bigoplus}_{i=0}^{n}CH^{i}(X)$$ be the Chow group of cycles of $X$ modulo
rational equivalence. Assume that $X$ is nonsingular over a field. Then

\medskip

a) There is a Chern class homomorphism
$$C_{t} : K_{0}(X)\rightarrow 1+\bigoplus_{i=1}^{n}CH^{i}(X)t^{i}$$
which is an additive to multiplicative group homomorphism. We write
$C_{t}(x) = 1+\sum_{i=1}^{n}C_{i}(x)t^{i}$ with $C_{i}(x)$ in $CH^{i}(X).$

\medskip

b)For nonnegative integer $k, {\Cal F}^{k}K_{0}(X)$ or simply ${\Cal F}^{k}(X)$ will denote the subgroup
of $K_{0}(X)$ generated by $[M]$, where $M$ runs through all coherent sheaves
on $X$ with $codimension(support M)$ atleast $k$. For such a coherent sheaf
$M$, $Cycle M$ will denote the codimention $r - cycle$ in the Chow group of
$X.$ (There will be no scope of confussion with notation $Cycle~x$ we
introduced in (6.1 c).) 

\medskip

c)We recall some of the properties of this Chern Class homomorphism
(see \cite{F}):
\newline
(1) if $x$ is in ${\Cal F}^{r}K_{0}(X)$ then 
$$C_{i}(x) = 0~~~	1\leq i < r	~~~  and$$
$$C_{r}(x) = (-1)^{r-1}(r-1)!Cycle(x)$$
\newline
(2)If $E$ is a locally free sheaf of finite rank over $X$ and there is a surjective
map from $E$ onto a locally complete intersection ideal sheaf $I$ of height $r$
then $$C_{r}(E) = (-1)^{r}Cycle({\Cal O}_{X}/I).$$

\endexample

	It is known that for a nonsingular variety $X$ over a field, the
$\gamma-filtration~F^{r}(X)$ of $K_{0}(X)$ is finer than the filtration
${\Cal F}^{r}(X)$ i.e. $F^{r}(X) \subseteq {\Cal F}^{r}(X)$. Following is an 
example of a nonsingular affine ring over a field $k$, for which these two 
filtrations indeed disagree.

\medskip

\example{Example 6.3} Following the notations in section 2,
 for a fixed positive
integer $n$ and a field $k$, let

$$
A_{n}=A_{n}(k)=\frac{k[S,T,U,V,X_{1}, \ldots ,X_{n},Y_{1}, \ldots ,Y_
{n}]}{(SU+TV-1,X_{1}Y_{1}+ \cdots +X_{n}Y_{n} - ST)}
$$

$$
B_{n}=B_{n}(k)=\frac{k[T,X_{1}, \ldots ,X_{n},Y_{1}, \ldots ,Y_{n}]}
{(X_{1}Y_{1} + \cdots +X_{n}Y_{n} -T(1+T))}
$$
Then for $X= Spec A_{n}$ or $Spec B_{n}$, 
$$ {\Cal F}^{r}(X) \approx {\Bbb Z}~~ for~~ 1 \leq r \leq n ~~ and $$ 
$${\Cal F}^{r}(X) = 0~~  for~~ n <r.$$
\newline
Further, $$ F^{n}(X)=(n-1)!{\Cal F}^{n}(X)
~~~ and~~~ F^{r}(X) = 0~~~ for~~~ n < r.$$
\endexample 

\demo{Proof} The computation of ${\Cal F}^{r}(X)$ is done exactly as in 
\cite{Sw1},(see (2.1) 
in case of $A_{n}$).Since $F^{r}(X)$ is contained in ${\Cal F}^{r}(X)$, $F^{r}(X)=0$ for $r > n$.

For definiteness, let $X$ be $Spec A_{n}$. So,
$\lambda_{n}= [A_{n}/(X_{1}, \ldots, X_{n},T)]$ is
the generator of ${\Cal F}^{r}(X)$ for $1 \leq r \leq n$.
Also $\lambda_{n}$ is the generator of $F^{1}(X)={\Cal F}^{1}(X)$.
By (2.6) in , there is 
a
projective $A-module~P$ of rank $n$ with $[P]-n = -\lambda_{n}$ so that $P$ 
maps onto the ideal $J=(X_{1},\ldots 
, X_{n-1})+I^{(n-1)!}$, where $I$ is the ideal $(X_{1},\ldots, X_{n},T)A_{n}$.
Hence $(n-1)!\lambda_{n} = -[A_{n}/J]$ ,(see (7.3)), is in $F^{n}(X)$.Also since
${\Cal F}^{n+1}(X)=0$, $\lambda_{n}^{2} = 0$.Hence for $1 \leq k$,$\gamma ^{k}$
acts as a group homomorphism on $F^{1}(X)$. 
So,$F^{n}(X)= {\Bbb Z}\gamma^{n}(\lambda_{n})$. As $F^{n+1}(X)=0$,
by (6.1), $\gamma_{n}((n-1)!\lambda_{n})= c_{n}((n-1)!\lambda_{n})=
(-1)^{n-1}(n-1)!^{2}\lambda_{n}$. Hence 
$F^{n}(X)={\Bbb Z}(n-1)!\lambda_{n}$. So the proof of (6.3) is complete.
\enddemo

\head 7. Some More Preliminaries\endhead

Following theorem (7.1) gives the Chern classes of the projective
module $P$ that we constructed in theorem (3.1).

\proclaim{Theorem 7.1} Under the set up and notations of theorem (3.1)
, we further have

$$c_{i}(P) = 0~~~ for~~~1 \leq i < r~~~ in~~~ \Gamma^{i}(X)\bigotimes {\Cal Q}$$,
$$c_{r} (P) = (-1)^{r}Cycle(A/J)~~~ in~~~ \Gamma^{r}(X).$$

If X is nonsingular over a field then
$$C_{i}(P) = 0~~~ for~~~1 \leq   i < r~~~ in~~~ CH^{i}(X)~~~ and$$
$$C_{r}(P) = (-1)^{r}Cycle(A/J)	~~~ in~~~ CH^{r}(X).$$
\endproclaim

\demo{Proof}  Comments about Chern
classes in $\Gamma(X)$ follow from (6.1), because $(r-1)!([P]-r)=-[A/J]$
is in $F^{r}(X)$. Similarly, since $[A/J_{0}]$ is in ${\Cal F}^{r} (X)$, the
comments about Chern classes in Chow group follows from (6.2).
\enddemo

\remark{Remark 7.2} For historical reasons we go back to the statement
of theorem 3.1. Let $J_{0}$ be an ideal in a noetherian 
commutative ring $A$ and
let $J_{0}=(f_{1},\ldots,f_{r})+J_{0}^{2}$.The part(i) of theorem (3.1) evolved
in two stages.First, Boratynski \cite{B} defined $J = (f_{1},\ldots,
f_{r-1}) + J_{0}^{(r-1)!}$ and proved that there is a projective $A-module~P$
of rank $r$ that maps onto $J$.Then, Murthy\cite{Mu2} added that if $J_{0}$ is a
locally complete intersection ideal of height $r$ then there is one such 
projective $A-module~P$ of rank $r$, with $[P]-r= -[A/J_{0}]$, that maps onto $J$.

We shall be much concerned with such ideals $J$ constructed, as above,
from ideals $J_{0}$. Following are some comments about such ideals.
\endremark

\proclaim{Natations and Facts 7.3} For an ideal $J$ 
in a Cohen-Macaulay ring $A$ with $J = (f_{1},
\ldots,f_{r}) + J^{2}$.We use the notation
 $$B(J)=B(J,f_{1},\ldots,f_{r}) = 
(f_{1},\ldots,f_{r-1})+J^{(r-1)!}.$$ Then, we have

$(1)~\sqrt{J}=\sqrt{B(J)},$

(2) $J$ is locally complete intersection ideal of height $r$ if and only if so
is $B(J)$.

(3) If $J$ is locally complete intersection ideal of height $r$ then 
$[A/B(J)] = (r-1)![A/J]$ in $K_{0}(X)$.
\endproclaim

\demo{Proof} The proof of (1) is obvious.To see (2), note that locally, $J_{0}$ is
generated by $f_{1}, \ldots ,f_{r}$ and $B(J_{0})$ is generated by $f_{1}, \ldots ,
f_{r-1}, f_{r}^{(r-1)!}.$ To prove (3), let $J_{k}=(f_{1}, \ldots , f_{(r-1)})
+J_{0}^{k}$ for positive integers $k$. Note that $0 \rightarrow J_{k}/J_{k+1}
 \rightarrow A/J_{k+1} \rightarrow A/J_{k} \rightarrow 0$ is exact and
$A/J_{0} \approx J_{k}/J_{k+1}.$ Now (3) follows by induction and hence the
proof of (7.3) is complete .\enddemo 

The following lemma describes such ideals $B(J_{0})$ very precisely.

\proclaim{Lemma 7.4} Let $A$ be a Cohen-Macaulay ring and $J$ be an ideal in $A$.
Then $J = B(J_{0}) = B(J_{0},f_{1},\ldots,f_{r})$ for some ideal $J_{0} =
(f_{1},\ldots,f_{r}) +J_{0}^{2}$ if and only if $J = 
(f_{1},\ldots , f_{r-1}, f_{r}^{(r-1)!}) + J^{2}$.
\endproclaim

\demo{Proof} To see the direct implication, let $J= (f_{1}, \ldots , f_{r-1})
+J_{0}^{(r-1)!}$ where $J_{0}= (f_{1}, \ldots , f_{r})+J_{0}^{2}.$ Then, it is
easy to check that $J=(f_{1}, \ldots , f_{r-1}, f_{r}^{(r-1)!})+J^{2}.$
Conversely, let $J=(f_{1},\ldots ,f_{r-1}, f_{r}^{(r-1)!})+J^{2}.$ By
Nakayama's lemma, there is an $s$ in $J$ such that
$$(1+s)J \subseteq (f_{1}, \ldots , f_{r-1}, f_{r}^{(r-1)!})~~ and~~~
J=(f_{1}, \ldots, f_{r-1}, f_{r}^{(r-1)!}, s).$$ Now we let $J_{0} = (f_{1},
\ldots ,f_{r}, s).$ It follows that $J=B(J_{0})$ and the proof of (7.4) is
complete.
\enddemo

The following lemma will be useful in the next section.

\proclaim{Lemma 7.5} Let $A$ be a Cohen-Macaulay ring of dimension $n$ and
let $I$ and $J$ be two locally complete intersection ideals of height $r$
with $I+J=A$.Then,
if $IJ=B({\Cal J})$ for some locally complete intersection ideal
${\Cal J}$ of height $r$ then $I=B(I_{0})$ for some locally complete intersection ideal $I_{0}$ of height $r$.
 Also if  $I= B(I_{0})~ and~J=B(J_{0})$ for locally complete intersection ideals 
$I_{0},J_{0}$ then $IJ=B({\Cal J})$ for some locally complete intersection ideal ${\Cal J}$ 
of height $r.$
\endproclaim

The proof is straightforward.

\remark{Remark} With careful formulation of the statements, the Cohen-Macaulay
condition  in (7.3), (7.4), (7.5) can be dropped.\endremark

\head 8.  Results on Chern classes\endhead

Our approach here is that if $Q$ is a projective module of rank $r$ over
a noetherian commutative ring $A$ then we try to  construct a projective $A-module~Q_{0}$ of rank $r-1$, so that the first $r-1$ Chern classes of
$Q$ are same as that of $Q_{0}$.Conversely, given a projective $A-module
~Q_{0}$ of rank $r-1$ and a locally complete intersection ideal $I$ of height
$r$, we attempt to construct a projective $A-module~Q$ of rank $r$ such that
the first $r-1$ Chern classes of $Q$ and $Q_{0}$ are same and the top Chern
class of $Q$ is $(-1)^{r}Cycle([A/I])$. Our first theorem(8.1) suggests that if
$r= n = dim X$, then for such a possibility to work, it is important that $[A/I]$ is 
divisible by $(n-1)!$.

\proclaim{Theorem 8.1} Let $A$ be noetherian commutative ring of dimension $n$ and
$X= Spec A$.Assume that $K_{0}(X)$ has no $(n-1)!$ torsion. Suppose that $Q$ is
a projective $A-module$ of rank $n$ and $Q_{0}$ is a projective $A-module$ of
rank $n-1$.Assume that the first $n-1$ Chern classes,in $\Gamma(X)$ (respectively,in $CH(X)$, if $X$ is nonsingular over a field), of $Q$ and $Q_{0}$ are same.
Then $\sum_{i=0}^{n}(-1)^{i}[\Lambda^{i}Q]$ is divisible by $(n-1)!$ in
$K_{0}(X)$. That means that if $Q$ maps onto a locally complete intersection ideal
$I$ of height $n$, then $[A/I]$ is divisible by $(n-1)!$ in $K_{0}(X)$.
\endproclaim

\demo{Proof} We can find a projective $A-module~P$ of rank $n$ such that $Q\bigoplus A^{n} \approx Q_{0} \bigoplus A \bigoplus P$.It follows that $c_{i}(P)= 0$
for $1 \leq i < n$ and $c_{n}(P) = c_{n}(Q)$ in $\Gamma(X)$. Write
$\rho = [P] - n$.

We claim that for $r = 0$ to $n-1$, $\beta_{r}\rho$ is in $F^{r+1} (X)$,
where $\beta_{r} = \Pi_{i=1}^{r-1} (i!)$.By Induction, assume that $\beta_{
r-1} \rho $ is in $F^{r}(X)$. Since $c_{r}(\beta_{r}\rho) = \beta_{r}c_{r} (\rho) = 0$,
also since $c_{r}(\beta_{r}\rho) = (-1)^{r-1}(r-1)!\beta_{r}\rho$ the
claim follows.

So, $\beta_{n-1}\rho$ is in $F^{n}(X)$.Hence 
$c_{n} (\beta_{n-1} \rho)
= (-1)^{(n-1)} (n-1)!\beta_{n-1} \rho$. Since $K_{0}(X)$ has no
 $\beta_{n-1}$ torsion, it follows that $c_{n}(\rho) = (-1)^{n-1}(n-1)!
\rho$.Since $c_{n}(\rho) = c_{n}(Q) = (-1)^{n}\sum_{i=0}^{n}(-1)^{i}
[\Lambda^{i}Q]$, the theorem follows.

We argue similarly when $X$ is nonsingular over a field and Chern classes
take values in the Chow gorup. In this case, we use 
(6.2 c). This completes the proof of (8.1).
\enddemo

Our next theorem(8.2) is a converse of (8.1).

\proclaim{Theorem 8.2} Let $A$ be a noetherian commutative ring of dimension $n$
and $X = Spec A$. Let $J$ be a locally complete intersection ideal of height
$r > 0$ with $J = (f_{1}, \ldots ,f_{r-1},f_{r}^{(r-1)!}) + J^{2}$ ( hence
$J = B(J_{0})$ for some locally complete intersection ideal $J_{0}$ of height
$r$). Let $Q$ be projective a $A-module$ of rank $r$ that maps onto $J$. Then
there is a projective $A-module~Q_{0}$ of rank $r-1$ such that,

(1)$~[Q_{0} \bigoplus A] = [Q] + [A/J_{0}]~~~ in~~~ K_{0}(X),$

(2)$$ c_{i}(Q_{0}) = c_{i}(Q)~~~ in~~~ \Gamma(X) \bigotimes {\Cal Q}~~~ for
~~~1 \leq i <r ~~~ and$$
\newline if $X$ is nonsingular over a field then
$$ C_{i}(Q_{0}) = C_{i}(Q)~~~ in~~~CH(X)~~~ for~~~ 1 \leq i < r.$$

(3) If $K_{0}(X)$ has no torsion (respectively, no $(n-1)!$ torsion,in case $X$
is nonsingular over a field) then such a $[Q_{0}]$ satisfying (2)
is unique in $K_{0}(X)$.
\endproclaim

\demo{Proof} By (3.1) with $I = A$, there is a surjective map 
$$\psi : Q \bigoplus A^{r} \rightarrow A \bigoplus P$$  where $P$ is a projective
$A- module$ of rank $r$ that maps onto $J$ and $[P] - r = - [A/J_{0}]$. We
let $Q_{0} = kernel\,( \psi)$. Then $Q_{0} \bigoplus A \bigoplus P \approx
Q \bigoplus A^{r}$. This settles (1). By (7.1), first $r-1$ Chern classes
of $P$, in $\Gamma(X)\bigotimes{\Cal Q}$, are zero. Hence it follows that
$c_{i}(Q) = c_{i}(Q_{0})$ for $1 \leq i < r$.
	In case $X$ is nonsingular, the argument runs similarly. So, the proof
of (2) of (8.2) is complete.

To prove (3), let $Q'$ be another projective $A-module$ of rank $r-1$
satisfying (2) and let $\rho = ([Q_{0}] - [Q'])$. Since the total Chern classes
of $Q_{0}$ and $Q'$ in $\Gamma(X)\bigotimes {\Cal Q}$(respectively,in $CH(X)$,
in case $X$ is nonsingular), are same,the total 
Chern class $c(\rho )=1$ in the respective groups.

For a positive integer $r$ let $\beta_{r}= \Pi_{i=1}^{r-1}(i!)$. By induction, as
in (8.1), it follows that $\beta_{n}\rho$ is in 
$F^{n+1}(X)\bigotimes {\Cal Q} = 0$ (respectively, in ${\Cal F}_{n+1}(X) = 0$).
Hence the proof 
of (8.2) is complete.
\enddemo

Following theorem(8.3) gives  a construction of projective modules with
certain given cycles as its total Chern class.

\proclaim{Theorem 8.3} Let $A$ be a Cohen- Macaulay ring
of dimension $n$ and $X = Spec A$. Let $r_{0}$ be an
integer with $2r_{0} \geq n$.

Let $Q_{0}$ be a projective $A-module$ of rank $r_{0}-1$, such that for
all locally complete intersection subschemes $Y$ of $X$
 with $codimension~Y \geq r_{0}$, the restriction $Q_{0}|Y$ of $Q_{0}$ to $Y$
 is trivial. Also
let $r$ be another integer with $r_{0} \leq r \leq n$ and for $k=r_{0}~ to~r$, let $I_{k}$ be locally 
complete intersection ideals of height $k$, with $I_{k}=(f_{1},\ldots,f_{k-1}
,f_{k}^{(k-1)!})+ I_{k}^{2}$ (hence $I_{k}=B(I_{k0})$, for some locally complete
intersection ideal $I_{k0}$ of height $k$).

Then there is a projective $A-module~Q_{r}$ of rank $r$ such that

(1) $Q_{r}$ maps onto $I_{r}$,

(2)
$$[Q_{r}] - r 
= ([Q_{0}]-(r_{0}-1)) + {[A/J_{r_{0}}] + \cdots + [A/J_{r}]},$$
where 
$J_{k}$ is a locally complete intersection ideal of height $k$ such that
\newline $(k-1)![A/J_{k}] = -[A/I_{k}]$ and further, $[P_{k}]-k = -[A/J_{k}]$ for some 
projective $A-module~P_{k}$ of rank $k$, for $r_{0} \leq k \leq r.$

(3)$$c_{k}(Q_{r}) =c_{k}(Q_{0})~~~ in~~~ \Gamma^{k}(X)\bigotimes {\Cal Q}~~~for
~~~1 \leq k < r_{o}$$
 $$c_{k}(Q_{r}) = (-1)^{k}Cycle([A/I_{k}]) ~~~in~~~ 
\Gamma^{k}(X)\bigotimes{\Cal Q}~~~ for~~~r_{o} \leq k \leq r.$$
\newline If $X$ is nonsingular over a field,then
$$C_{k}(Q_{r})=C_{k}(Q_{0})~~~in~~~ CH^{k}(X)~~~for~~~1 \leq k <~r_{0}~~~ and$$
$$C_{k}(Q_{r})= (-1)^{k}Cycle(A/I_{k})~~~in~~~CH^{k}(X)~~~for
~~~r_{0} \leq k \leq r.$$
\endproclaim

{\bf Caution.} In the statement of (8.3) the generators $f_{1} , \ldots ,f_{k-1}, f_{k}^{(k-1)!}$ of $I_{k}/I_{k}^{2}$ depend on $k$.

\remark{Remark 8.4} A free $A-module~Q_{0}$ of rank $r_{0}-1$ will satisfy the
hypothesis of (8.3). If $r_{0} = n-1$, then any projective $A-module~Q_{0}$
of rank $r_{0}-1$ with trivial determinant will also satisfy the hypothesis
of (8.3). 
\endremark

The proof of (8.3) follows, by induction, from the following proposition(8.5).

\proclaim{Proposition 8.5} Let $A$ be a Cohen-Macaulay ring 
of dimension $n$ and $X = Spec A$ and let $r$ be
a positive integer with $2r \geq n$ and $r \leq n$. Let $Q_{0}$ be a projective $A-module$ of
rank $r-1$ such that for any locally complete intersection closed subscheme
$Y$ of codimension atleast $r$, the restriction $Q_{0}|Y$ of $Q_{0}$ to $Y$ is
trivial. Also let $I$ be a locally complete intersection ideal of height $r$
with $I = (f_{1},\ldots,f_{r-1},f_{r}^{(r-1)!}) + I^{2}$ (hence $I= B(I_{0})$
for some locally complete intersection ideal $I_{0}$ of height $r$).

Then there is a projective $A-module~ Q$ of rank $r$ such that

(1) $Q$ maps onto $I$,

(2)$$[Q] - r = ([Q_{0}] - (r-1)) + [A/J_{0}],$$ where $J_{0}$ is a locally
complete intersection ideal of height $r$ such that $(r-1)![A/J_{0}] = -[A/I]$
and further there is a projective $A-module~P$ of rank $r$ such that $[P] - r =
-[A/J_{0}]$.

(3)$$c_{k}(Q)=	c_{k}(Q_{0})~~~in~~~ \Gamma^{k}(X)\bigotimes{\Cal Q}~~~ for
~~~ 1 \leq k < r,~~~and$$

$$ c_{k}(Q)= (-1)^{r}Cycle([A/I])~~~~in~~~ \Gamma^{k}(X)\bigotimes{\Cal Q}~~~
for~~~ k = r.$$

If $X$ is nonsingular over a field, then

$$ C_{k}(Q) =	C_{k}(Q_{0})~~~in~~~ CH^{k}(X)~~~for~~~ 1 \leq k < r$$ 
and
$$ C_{k}(Q)=(-1)^{r}Cycle(A/I)~~~in~~~CH^{k}(X)~~~for~~~ k =r.$$

(4) For any locally complete intersection closed subscheme $Y$ of $X$ of codimension atleast $r+1$, the restriction $Q|Y$ is trivial.
\endproclaim

Before we prove (8.5), we state the following proposition from \cite{CM}.

\proclaim{Proposition 8.6} Let $A$ be a noetherian commutative ring and $J$ be a 
locally complete intersection ideal of height $r$ with $J/J^{2}$ free. Suppose
I is an ideal with $dim A/I < r$. If $\pi :K_{0}(A) \rightarrow K_{0}(A/I)$
is the natural map then $\pi([A/J]) = 0$.
\endproclaim

The proof is done by finding a locally complete intersection ideal $J'$
of height $r$ such that $J'+I = A = J'+J$ and $J \cap J'$ is complete 
intersection. Now it follows that $\pi([A/J]) = -\pi([A/J']) = 0$.

\demo{Proof of (8.5)} We have $$I = (f_{1},\ldots,f_{r-1},f_{r}^{(r-1)!}) + I^{2}
=B(I_{0}) = (f_{1}, \ldots ,f_{r-1}) + I_{0}^{(r-1)!},$$ where $I_{0}$ is a locally
complete intersection ideal of height $r$ with $I_{0} = (f_{1},\ldots,f_{r})+
I_{0}^{2}$. We can also assume that $f_{1},\ldots,f_{r}$ is a regular sequence.
By hypothesis $Q_{0}/IQ_{0}$ is free of rank $r-1$. Let $e_{1},\ldots,e_{r-1}$
be elements in $Q_{0}$ whose images forms a basis of $Q_{0}/IQ_{0}$.So, there is
a map $\phi_{0}: Q_{0}\rightarrow I$ such that $\phi_{0}(e_{i})- f_{i}$ is in 
$I^{2}$ for $ i = 1$ to $r-1$.

So, $(\phi_{0}(Q_{0}), f_{r}^{(r-1)!})+I^{2}=I$. By Nakayama's lemma there is an 
$s$ in $I$ such that 
$$(1+s)I \subseteq (\phi_{0}(Q_{0}),f_{r}^{(r-1)!})~~~and~~~ 
I=
(\phi_{0}(Q_{0}),f_{r}^{(r-1)!},s).$$

Let $Q_{0}^{*}$ be the dual of $Q_{0}$. Then $(\phi_{0},s^{2})$ is 
basic in $Q_{0}^{*}\bigoplus A$ on the set ${\Cal P} = \{\wp~in~ Spec A\, :\, f_{r}~is
~in~ \wp~ and~height~(\wp) \leq r-1\}$. There is a generalised dimension
function $d\, :\, {\Cal P} \rightarrow \{0,1,2,\ldots\}$ so that $d(\wp) \leq r-2$
for all $\wp$ in ${\Cal P}$ (see \cite{P}). Since rank $Q_{0}^{*} = r-1 > d(\wp
)$ for all $\wp$ in ${\Cal P}$, there is an $h$ in $Q_{0}^{*}$ such that $\phi = \phi_{0}+s^{2}h$ is basic in $Q_{0}^{*}$ on ${\Cal P}$.

Write ${\Cal I} = (\phi(Q_{0}),f_{r}^{(r-1)!})$. It follows that (1)
${\Cal I}$ is a locally complete intersection ideal of height $r$, (2) $[A/{\Cal I}]=0,$
(3)  $ {\Cal I}+I^{2}=I$ (4) ${\Cal I} =
(g_{1},\ldots,g_{r-1},f_{r}^{(r-1)!})+{\Cal I}^{2}$ for some $g_{1},\ldots,g_{r-1}
$ in ${\Cal I}$.

To see (1), note that ${\Cal I}$ is locally $r$ generated and also since $\phi$ is basic in $Q_{0}^{*}$ on $\Cal P$, $\Cal I$ has height $r$.Now
since $A$ is Cohen-Macaulay, $\Cal I$ is a locally complete intersection ideal of height $r$. Since 
$$ 0 \rightarrow A/\phi(Q_{0}) @>f_{r}^{(r-1)!}>> A/\phi(Q_{0}) \rightarrow
 A/{\Cal I} \rightarrow 0$$
is exact, (2) follows. Since $\phi=\phi_{0}+s^{2}h$,(3) follows. By
hypothesis $Q_{0}/{\Cal I}Q_{0}$ is free of rank $r-1$ and hence (4) follows.

Because of (4), ${\Cal I} = B({\Cal I}_{0})$ for some locally complete intersection ideal ${\Cal I}_{0}$ of height $r$. From (3) it follows that ${\Cal I}=J\cap I$ for
some locally complete intersection ideal $J$ of height $r$ and $I+J=A$.
Since  ${\Cal I}=B({\Cal I}_{0})$, by (7.5), $J=B(J_{0})$ for some locally complete intersection
 ideal $J_{0}$ of height $r$.

Let $\phi : Q_{0}\bigoplus A \rightarrow {\Cal I}$ be the surjective map 
$(\phi,f_{r}^{(r-1)!})$. We can apply theorem(3.1) and (7.1).
 There is a surjective map 
$\psi:Q_{0}\bigoplus A^{r+1} \rightarrow P \bigoplus I$, where $P$ is a projective $A-module$ 
of rank $r$ that maps onto $J$ and $[P]-r =-[A/J_{0}]$. Also
$$c_{k}(P) = 0~~~for~~~	1 \leq k < r~~~in~~~ \Gamma^{k}(X)\bigotimes{\Cal Q},$$
$$c_{r}(P) = (-1)^{r}Cycle(A/J)~~~	in~~~ \Gamma^{r}(X).$$
\newline If $X$ is nonsingular over a field then
$$C_{k}(P) = 0~~~ for ~~~	1 \leq k < r~~~	in~~~ CH^{k}(X)~~~ and~~~$$
$$C_{r}(P) = (-1)^{r}Cycle(A/J)~~~	in~~~ CH^{r}(X).$$
Now $Q = \psi^{-1}(I)$ will satisfy the assertions of the theorem.Clearly,
$Q$ maps onto $I$ and (1) is satisfied. Note that $Q \bigoplus P \approx Q_{0}
\bigoplus A^{r+1}$ and hence
$$[Q]-r = ([Q_{0}]-(r-1))-([P]-r) = ([Q_{0}]-(r-1))+[A/J_{0}].$$
Also $(r-1)![A/J_{0}]=[A/J]= -[A/I]$, since $[A/{\Cal I}] = 0$.This establishes (2).	

Again since $Q \bigoplus P \approx Q_{0} \bigoplus A^{r+1}$ and since the Chern 
classes of $P$ are given as above, (3) follows.

To see (4), let $Y$ be a locally complete intersection subscheme of $X$
with codimension at least $r+1$.Let $\pi:K_{0}(X)\rightarrow K_{0}(Y)$ be the
restriction map. Then $\pi ([Q]-r) = \pi ([Q_{0}]-(r-1))+\pi([A/J_{0}])=0$ by
(8.6).Hence the restriction $Q|Y$ is stably free.Since $r > dim Y$, by cancellation theorem of Bass(see \cite{EE}), $Q|Y$ is free. This completes the proof of  
(8.5).
\enddemo 

Before we go into some of the applications let us recall(1.5, 1.6) that
for a smooth affine variety $X=Spec A$ of dimension $n$ over a field, ${\Cal F}^{n}(X) = 
F_{0}K_{0}(X) = \{[A/I]~in~ K_{0}(X): I~ is~a~locally~complete~intersection
~ideal~of~height~n\}.$

\medskip

Following is an important corollary to theorem(8.3).

\proclaim{Corollary 8.7} Suppose $X=SpecA$ is a smooth affine variety of dimension
$n$ over a field. Assume that $CH^{n}(X)$ is divisible by $(n-1)!$. Let $Q_{0}$
be projective $A-module$ of rank $n-1$ and $x_{n}$ is a cycle in $CH^{n}(X)$.
Then there is a projective $A-module~Q$ of rank $n$ such that
$$C_{i}(Q)=C_{i}(Q_{0}) ~~~for~~~ 1\leq i < n~~~ and~~~ C_{n}(Q) =x_{n}~~~ in~~~ CH^{n}(X).$$
	Conversely, if $Q$ is a projective $A-module$ of rank $n$, then there
is a projective $A-module~Q'$ of rank $n$ such that 
$$C_{i}(Q)=C_{i}(Q')~~~for~~~ 1\leq i < n~~~and~~~ C_{n}(Q') = 0 ~~~in~~~ CH^{n}(X).$$
\endproclaim

\demo{Proof}
Since
the Chern class map $C_{n}:{\Cal F}^{n}(X) \rightarrow CH^{n}(X)$ sends $[A/I]$
to 
\newline $(-1)^{(n-1)}(n-1)!Cycle(A/I)$ (see \cite{F}), this map is surjective.
Since ${\Cal F}^{n}(X)=F_{0}K_{0}(X)$, there is a locally complete intersection
ideal $I_{0}$ of height $n$ such that $C_{n}(A/I_{0}) = -x_{n}$. By theorem(8.3)
with $I=B(I_{0})$, there is a projective $A-module~Q$ of rank $n$ such that
$[Q]-n = ([Q_{0}]-(n-1))+[A/J]$ where $J$ is a locally complete intersection
ideal of height $n$ with $(n-1)![A/J]=-[A/I]=-(n-1)![A/I_{0}]$.Hence 
$$C_{i}(Q)
=C_{i}(Q_{0}) ~~~for~~~ 1 \leq i < n ~~~and$$ 
$$C_{n}(Q)=C_{n}([A/J]) = (-1)^{n-1}(n-1)!Cycle([A/J)=$$ 
$$(-1)^{n-1}Cycle((n-1)![A/J])=(-1)^{n-1}Cycle(-(n-1)![A/I_{0}]) = x_{n}.$$This establishes the direct implication.

To see the converse, note that, as above, there is a projective 
$A$-module $ P$ such that $C_{i}(P) = 0 ~for~ 1 \leq i < n$ and $C_{n}(P)=-C_{n}(Q)$.
 Now $Q\bigoplus P \approx Q'\bigoplus A^{n}$ for some projective $A-module~Q'$
of rank $n$. It is obvious that $Q'$ satisfies the assertions. This completes
the proof of (8.7).
\enddemo

Following theorem of Murthy(\cite{Mu2}) follows from (8.7).

\proclaim{Theorem 8.8(Murthy)} Let $X = Spec A$ be a smooth affine variety of
dimension $n$ over an algebraically closed field $k$. Let $x_{i}$ be cycles 
in $CH^{i}(X)$ for $1 \leq i \leq n$. Then there is a projective $A-module
~Q_{0}$ of rank $n-1$ with the total Chern class $C(Q_{0}) = 1 + x_{1} +\cdots+
x_{n-1}$ in $CH(X)$ if and only if there is a projective $A-module~Q$ of rank
$n$ with the toal Chern class $C(Q) = 1 + x_{1} + \cdots + x_{n-1} + x_{n}.$
\endproclaim

\demo{Proof} In this case $CH^{n}(X)$ is divisible (see \cite{Le,Sr, Mu2}).
So the direct implication is immediate from (8.7).

To see the converse, let $Q'$ be as in (8.7). Since $C_{n}(Q')=0$,it
follows from the theorem of Murthy(\cite{Mu2}) that $Q'\approx Q_{0}\bigoplus A$ 
for some projective $A-module~Q_{0}$ of rank $n-1$.
It is obvious that $C(Q_{0})=1+x_{1}+\cdots+x_{n-1}$. So the proof of (8.8) is
complete.
\enddemo

Following is an alternative proof of the theorem of Mohan Kumar and
Murthy (\cite{MM}).

\proclaim{Theorem 8.9 (\cite{MM})}  Let $X=Spec A$ be a smooth affine three fold over an
algebraically closed field and let $x_{i}$ be cycles in  $CH^ {i}(X)$ for $1 \leq i \leq 3.$
Then

(1) There is projective $A-module~Q_{3}$ of rank $3$ with total Chern class
$C(Q_{3}) = 1+x_{1}+x_{2}+x_{3},$

(2)there is a projective $A-module~Q_{2}$ of rank $2$ with total Chern calss
$C(Q_{2})=1=x_{1}+x_{2}$.\endproclaim

\demo{Proof} Because of (8.8), we need to prove (1) only.Let $L$ be a line
bundle on $X$ with $C_{1}(L)=x_{1}$.We claim that there is a projective 
$A-module~P$ of rank $3$ so that $C_{1}(P) = 0$ and $C_{2}(P) = x_{2}$.Let $x_{2}=
(y_{1}+\cdots +y_{r}) - (y_{r+1}+\cdots +y_{s})$ where
$y_{i}$ is the cycle of $A/I_{i}$ for prime ideals $I_{i}$ of height $2$ for $1 \leq i \leq s$.

For $1 \leq i \leq s$ there is an exact sequence
$$0\rightarrow P_{i}\bigoplus G_{i} \rightarrow F_{i} \rightarrow A \rightarrow A/I_{i} \rightarrow 0$$
where $F_{i}$ and $G_{i}$ are free modules and $P_{i}$ are projective $A-modules
$ of rank $3$.Since the total Chern classes $C(P_{i})=C(A/I_{i})$, it follows
that $C_{1}(P_{i})
=C_{1}(A/I_{i})= 0$ and $C_{2}(P_{i})= C_{2}(A/I_{i})= - Cycle (A/I_{i}) =-y_{i} 
$. There are free modules R and S and a projective $A-module~P$ of rank $3$ such
that $P_{1}\bigoplus\cdots P_{r}\bigoplus R \bigoplus P \approx P_{r+1}\bigoplus
\cdots \bigoplus P_{s}\bigoplus S$. It follows that $C_{1}(P) = 0$ and $C_{2}(P)
= (y_{1}+\cdots+y_{r})-(y_{r+1}+\cdots+y_{s})= x_{2}$. This establishes the
claim.

Let $P\bigoplus L \approx P'\bigoplus A$. Then $C_{1}(P')=C_{1}(L)=x_{1}$ and
$C_{2}(P')=x_{2}$ and let $C_{3}(P')=z$ for some $z$ in $CH^{3}(X)$. Again,by
(8.7) there is a projective $A-module~Q'$ of rank $3$ such that the total
Chern class $C(Q')=1+(x_{3}-z)$. There is a projective $A-module~Q_{3}$ of
rank $3$ such that $Q'\bigoplus P' \approx Q_{3}\bigoplus A^{3}$.We have,
$C(Q')C(P')=C(Q_{3})$. So the proof of (8.9) is complete.\enddemo

The same proof of (8.9) yeilds the following stronger theorem (8.10).

\proclaim{Theorem 8.10} Let $X = Spec A$ be a smooth affine three 
fold over any field $k$ such that $CH^{3}(X)$ is divisible by two.
Given $x_{i}$ in $CH^{i}(X)$ for $1 \leq i \leq 3$, there is a 
Projective $A-module~ Q$ of rank 3 such that the total chern class
$C(Q) = 1 + x_{1} + x_{2} + x_{3}.$\endproclaim

\remark{Remark.} For examples of smooth three folds that satisfy
the hypothesis of (8.10) see \cite {Mk2}.

\Refs

\refi B M. Boratynski, \it A note on set-theoretic complete intersection
ideals, \rm J. Algebra  {\bf 54}(1978).

\refi CF L. Claborn and R. Fossum, \it Generalization of the notion of class group, \rm 
Ill. J. of Math. {\bf 12}(1968), 228-253.

\refi CM Fernando Cukierman and Satya Mandal,
\it Study of vector bundles by restriction,
\rm to appear in Comm. Algebra.

\refi EE D. Eisenbud and E. G. Evans,
\it Generating modules efficiently: theorems from algebraic K-theory,
\rm J. Algebra {\bf 49}(1977) 276-303.

\refi F William Fulton,
\it Intersection Theory,
\rm Springer-Verlag(1984).

\refi FL William Fulton and Serge Lang,
\it Riemann-Roch Algebra,
\rm Springer-Verlag(1985)

\refi SGA6 A. Grothendieck et al, \it  
Theorie des Intersections et Theoreme de Riemann-Roch, 
\rm LNM {\bf 225}, Springer-Verlag (1971).

\refi Hu Dale Husemoller, \it Fibre Bundle, \rm GTM {\bf 20}, Springer-Verlag (1966).

\refi K Ernst Kunz, \it Kahler Differentials, \rm Vieweg (1986).

\refi J J. P. Jounanlou, \it Quelques calculs en K-theory des schemas, \rm Algebraic
K-Theory I, LNM {\bf 341} \rm Springer-Verlag, Berlin (1973), 317-334.

\refi Le Marc Levine, \it Zero cycles and K-theory on singular varieties, \rm Algebraic
Geometry, Bowdoin, Proc. Symp. in Pure Math {\bf 46}(1987), 451-462.

\refi Mk1 N. Mohan Kumar, \it Some theorems on generation of ideals in affine algebras, \rm
Comm. Math. Helv. {\bf 59}(1984), 243-252.

\refi Mk2 N. Mohan Kumar, \it Stably free modules, \rm Amer. J. of Math., {\bf 107}(1977),
1439-1443.

\refi MM N. Mohan Kumar and M. P. Murthy, \it Algebraic cycles and vector bundles over affine
three-folds, \rm Annals of Math. {\bf 116}(1982), 579-591.

\refi Mn Y. I. Manin,
\it Lectures on K-Functors in algebraic geometry,
\rm Russ. Math. Surveys {\bf 24},no 5(1969), 1-89.

\refi Mu1 M. P. Murthy, \it Zero-cycles, splitting of projective modules and number of
generators of modules, \rm Bull. Amer. Math. Soc. {\bf 19}(1988), 315-317.

\refi Mu2 M. P. Murthy, \it Zero cycles and projective modules,
 \rm Ann. of Math. {\bf 140}(1994), 405-434.

\refi Mu3 M. P. Murthy, \it Vector bundles over affine surfaces birationally equivalent to
ruled surfaces, \rm Ann. of Math. {\bf 89}(1969), 242-253.

\refi P B. R. Plumstead, \it The conjectures of Eisenbud and Evans, 
\rm Amer. J. of Math. {\bf 105} (1983), 1417-1433.

\refi Q Daniel Quillen, \it Higher Algebraic K-Theory, \rm Algebraic K-Theory I, LNM {\bf 341},
Springer-Verlag (1973), 85-147.

\refi Sr V. Srinivas, \it Torsion 0-cycles on affine varieties in characteristic p, \rm J. of
Alg. {\bf 120}(1989), 428-432.

\refi S A. A. Suslin, \it On stably free modules, \rm Math. USSR Stornik, {\bf 40}(1977).

\refi Sw1 R. G. Swan, \it Vector bundles, projective modules and the K-theory of spheres, \rm
Proc. of John Moore Conference, Ann. of Math. Study {\bf 113}(1987), 432-522.

\refi Sw2 R. G. Swan, \it K-Theory of quartic hypersurfaces, \rm Ann. of Math {\bf 122}(1985),
113-153.
\endRefs

\end